\documentclass[10pt]{article}%
\usepackage{amsfonts}
\usepackage{amssymb}
\usepackage{amsthm}
\usepackage{amsmath}
\usepackage{graphicx}%
\theoremstyle{plain}
\newtheorem{thm}{Theorem}[section]
\newtheorem{lem}[thm]{Lemma}
\newtheorem{cor}[thm]{Corollary}
\newtheorem{prop}[thm]{Proposition}

\theoremstyle{definition}
\newtheorem{dfn}[thm]{Definition}

\newtheorem{ex}{Example}

\begin{document}

\title{Index Sets of Computable Structures }
\author{Wesley Calvert\thanks{Calvert acknowledges the support of the NSF grants
DMS-0139626 and DMS-0353748, Harizanov of the NSF grant DMS-0502499, and
Miller of the NSF grant DMS-0353748. Harizanov also gratefully acknowledges
the support of the Columbian Research Fellowship of the George Washington
University.}\\Department of Mathematics \& Statistics \\Murray State University\\Murray, Kentucky 42071\\wesley.calvert@murraystate.edu
\and Valentina S. Harizanov\\Department of Mathematics \\George Washington University\\harizanv@gwu.edu
\and Julia F. Knight\\Department of Mathematics \\University of Notre Dame\\knight.1@nd.edu
\and Sara Miller\\Department of Mathematics \\University of Notre Dame\\smiller9@nd.edu}
\maketitle

\begin{abstract}
The \emph{index set} of a computable structure $\mathcal{A}$ is the set of
indices for computable copies of $\mathcal{A}$. We determine the complexity of
the index sets of various mathematically interesting structures, including
arbitrary finite structures, $\mathbb{Q}$-vector spaces, Archimedean real
closed ordered fields, reduced Abelian $p$-groups of length less than
$\omega^{2}$, and models of the original Ehrenfeucht theory. The index sets
for these structures all turn out to be $m$-complete $\Pi_{n}^{0}$,
$d$-$\Sigma_{n}^{0}$, or $\Sigma_{n}^{0}$, for various $n$. In each case, the
calculation involves finding an \textquotedblleft optimal\textquotedblright%
\ sentence (i.e., one of simplest form) that describes the structure. The form
of the sentence (computable $\Pi_{n}$, $d$-$\Sigma_{n}$, or $\Sigma_{n}$)
yields a bound on the complexity of the index set. When we show $m$%
-completeness of the index set, we know that the sentence is optimal. For some
structures, the first sentence that comes to mind is not optimal, and another
sentence of simpler form is shown to serve the purpose. For some of the
groups, this involves Ramsey theory.

\end{abstract}

\section{Introduction}

One of the goals of computable structure theory is to study the relationship
between algebraic and algorithmic properties of structures. Our languages are
computable, and our structures have universes contained in $\omega$, which we
think of as computable sets of constants. If $\varphi$ is a formula, we write
$+\varphi$ for $\varphi$, and $-\varphi$ for $\lnot\varphi$. In measuring
complexity, we identify a structure $\mathcal{A}$ with its atomic diagram,
$D(\mathcal{A)}$. In particular, $\mathcal{A}$ is \emph{computable} if
$D(\mathcal{A})$ is computable.

For a computable structure $\mathcal{A}$, an \emph{index} is a number $a$ such
that $\varphi_{a}=\chi_{D(\mathcal{A})}$, where $(\varphi_{a})_{a\in\omega}$
is a computable enumeration of all unary partial computable functions. The
\emph{index set} for $\mathcal{A}$ is the set $I(\mathcal{A})$ of all indices
for computable (isomorphic) copies of $\mathcal{A}$. For a class $K$ of
structures, closed under isomorphism, the \emph{index set} is the set $I(K)$
of all indices for computable members of $K$. There is quite a lot of work on
index sets \cite{LS}, \cite{CMS}, \cite{C2}, \cite{C1}, \cite{CCKM},
\cite{G-K}, \cite{White}, \cite{WhiteD}, \cite{dobritsa}, etc. Our work is
very much in the spirit of Louise Hay, and Hay together with Doug Miller (see
\cite{DM}).

In this paper, we present evidence for the following thesis:

\begin{quotation}
\noindent For a given computable structure $\mathcal{A}$, to calculate the
precise complexity of $I(\mathcal{A})$, we need a good description of
$\mathcal{A}$, and once we have an \textquotedblleft optimal\textquotedblright%
\ description, the complexity of $I(\mathcal{A})$ will match that of the description.
\end{quotation}

\noindent Our evidence for the thesis consists of calculations for computable
structures of several familiar kinds: finite structures, $\mathbb{Q}$-vector
spaces, Archimedean ordered fields---the ones we consider are real closed or
purely transcendental extensions of $\mathbb{Q}$, reduced Abelian $p$-groups
of length less than $\omega^{2}$, and models of the original Ehrenfeucht theory.

We should say what qualifies as a \textquotedblleft
description\textquotedblright\ of a structure, and how we measure the
complexity. The Scott Isomorphism Theorem says that for any countable
structure $\mathcal{A}$, there is a sentence of $L_{\omega_{1}\omega}$ whose
countable models are exactly the isomorphic copies of $\mathcal{A}$ (see
\cite{keisler}). Such a sentence is called a \emph{Scott sentence} for
$\mathcal{A}$. A Scott sentence for $\mathcal{A}$ certainly describes
$\mathcal{A}$.

There is earlier work \cite{DM}, \cite{AM} investigating subsets of the Polish
space of structures with universe $\omega$ for a given countable relational
language. Concerning the possible complexity (in the noneffective Borel
hierarchy) of the set of copies of a given structure, it is shown in \cite{DM}
that if the set is $\mathbf{\Delta_{\alpha+1}^{0}}$, then it is $d$%
-$\mathbf{\Sigma_{\alpha}^{0}}$. In \cite{AM} it is shown that the set cannot
be properly $\mathbf{\Sigma_{2}^{0}}$. There are also examples illustrating
other possibilities.

Most of the structures we consider follow one of two patterns. Either there is
a computable $\Pi_{n}$ Scott sentence, and the index set is $m$-complete
$\Pi_{n}^{0}$, or else there is a Scott sentence which is computable
\textquotedblleft$d$-$\Sigma_{n}$\textquotedblright\ (the conjunction of a
computable $\Sigma_{n}$ sentence and a computable $\Pi_{n}$ sentence), and the
index set is $m$-complete $d$-$\Sigma_{n}^{0}$. For example, a computable
reduced Abelian $p$-group of length $\omega$ has a computable $\Pi_{3}$ Scott
sentence, and the index set is $m$-complete $\Pi_{3}^{0}$. A $\mathbb{Q}%
$-vector space of finite dimension at least $2$ has a Scott sentence that is
computable $d$-$\Sigma_{2}$, and the index set is $m$-complete $d$-$\Sigma
_{2}^{0}$.

The \textquotedblleft middle model\textquotedblright\ of the original
Ehrenfeucht theory illustrates a further pattern. There is a computable
$\Sigma_{3}$ Scott sentence, and the index set is $m$-complete $\Sigma_{3}%
^{0}$. Often, the first Scott sentence that comes to mind is not optimal. In
some cases, in particular, for some of the groups, it requires effort to show
that a certain sentence of a simpler form actually is a Scott sentence.

For some structures, we obtain more meaningful results by locating the given
computable structure $\mathcal{A}$ within some natural class $K$. We say how
to describe $\mathcal{A}$ within $K$, and also how to calculate the complexity
of $I(\mathcal{A})$ within $K$.

\begin{dfn}
A sentence $\varphi$ is a Scott sentence for $\mathcal{A}$ \emph{within $K$}
if the countable models of $\varphi$ in $K$ are exactly the isomorphic copies
of $\mathcal{A}$.
\end{dfn}

The following definitions were already used in \cite{C2}.

\begin{dfn}
Let $\Gamma$ be a complexity class (e.g., $\Pi^{0}_{3}$).

\begin{enumerate}
\item $I(\mathcal{A})$ is $\Gamma$ \emph{within $K$} if $I(\mathcal{A}) =
R\cap I(K)$ for some $R\in\Gamma$.

\item $I(\mathcal{A})$ is $m$-\emph{complete} $\Gamma$ \emph{within $K$} if
$I(\mathcal{A})$ is $\Gamma$ within $K$ and for any $S\in\Gamma$, there is a
computable function $f:\omega\rightarrow I(K)$ such that
\[
n\in S\text{ iff}\ f(n)\in I(\mathcal{A})\text{;}%
\]
that is, there is a uniformly computable sequence $(\mathcal{C}_{n}%
)_{n\in\omega}$ for which
\[
n\in S\text{ iff}\ \mathcal{C}_{n}\cong\mathcal{A}\text{.}%
\]

\end{enumerate}
\end{dfn}

\begin{ex}
Let $\mathcal{A}$ be a field with $3$ elements, and let $K$ be the class of
finite prime fields. There is a Scott sentence for $\mathcal{A}$ within $K$
saying $1+1+1=0$. The index set for $\mathcal{A}$ is computable \emph{within}
$K$.
\end{ex}

The example above is an exception. In most of the examples we consider, even
when we locate our structure within a class $K$, the optimal description is a
true Scott sentence, but the context helps us calculate the complexity of the
index set in a meaningful way.

\begin{ex}
Let $\mathcal{A}$ be a linear ordering of size $3$, and let $K$ be the class
of linear orderings. There is a computable $d$-$\Sigma_{1}$ Scott sentence
saying that there are at least $3$ elements ordered by the relation, and not
more. We will show that the index set for $\mathcal{A}$ is $m$-complete
$d$-c.e.\ \emph{within $K$}.
\end{ex}

Here we mention some related work. The proof of the Scott Isomorphism Theorem
leads to an assignment of ordinals to countable structures. By a result of
Nadel \cite{Nadel}, for any hyperarithmetical structure, there is a computable
infinitary Scott sentence iff the Scott rank is computable. Several different
definitions of Scott rank are used. Since we are more interested in Scott
sentences, we shall not give any of them.

Work on index sets for particular computable structures is related to work on
isomorphism problems for classes of computable structures \cite{C1},
\cite{C2}, \cite{G-K}. The \emph{isomorphism problem} for a class $K$ is the
set $E(K)$ consisting of pairs $(a,b)$ of indices for computable members of
$K$ that are isomorphic. It is often the case that for the classes $K$ for
which the complexity of the isomorphism problem is known, there is a single
computable $\mathcal{A}\in K$ such that the index set for $\mathcal{A}$ has
the same complexity as $E(K)$. Results on index sets are useful in other
contexts as well. In \cite{CCHM}, they are used in connection with $\Delta
_{2}^{0}$ categoricity of computable structures.

The results on finite structures are in Section 2, and those on vector spaces
are in Section 3. In Section 4, we consider Archimedean real closed ordered
fields. The results on reduced Abelian $p$-groups are in Section 5, and the
results on models of the original Ehrenfeucht theory are in Section 6.

\section{Finite structures}

\label{fssec}

Finite structures are the easiest to describe. It is perhaps surprising that
there should be any variation in complexity of index sets for different finite
structures, and, indeed, there is almost none. In the following theorem, we
break with convention by allowing a structure to be empty.

\begin{thm}
\label{finite}

Let $L$ be a finite relational language. Let $K$ be the class of finite
$L$-structures, and let $\mathcal{A}\in K$.

\begin{enumerate}
\item If $\mathcal{A}$ is empty, then $I(\mathcal{A})$ is $m$-complete
$\Pi^{0}_{1}$ within $K$.

\item If $\mathcal{A}$ has size $n\geq1$, then $I(\mathcal{A})$ is
$m$-complete $d$-c.e.\ within $K$.
\end{enumerate}
\end{thm}

\proof

For 1, first note that $\mathcal{A}$ has a finitary $\Pi_{1}$ Scott sentence
saying that there is no element. From this, it is clear that $I(\mathcal{A})$
is $\Pi_{1}^{0}$ within $K$. For completeness, let $\mathcal{B}$ be an
$L$-structure with just one element. For an arbitrary $\Pi_{1}^{0}$ set $S$,
we can produce a uniformly computable sequence $(\mathcal{A}_{n})_{n\in\omega
}$ such that
\[
\mathcal{A}_{n}\cong\left\{
\begin{array}
[c]{ll}%
\mathcal{A} & \ \text{if}\ n\in S\text{,}\\
\mathcal{B} & \ \text{if}\ n\notin S\text{.}%
\end{array}
\right.
\]

For 2, we have a finitary existential sentence $\varphi$ stating that there is
a substructure isomorphic to $\mathcal{A}$, and another finitary existential
sentence $\psi$ stating that there are at least $n+1$ elements. Then
$\varphi\ \&\ \lnot\psi$ is a Scott sentence for $\mathcal{A}$. It follows
that $I(\mathcal{A})$ is $d$-c.e.\ within $K$. For completeness, let $S =
S_{1} - S_{2}$, where $S_{1}$ and $S_{2}$ are c.e. We have the usual finite
approximations $S_{1,s}$, $S_{2,s}$.

Let $\mathcal{A}^{-}$ be a proper substructure of $\mathcal{A}$, and let
$\mathcal{A}^{+}$ be a finite proper superstructure of $\mathcal{A}$. We will
build a uniformly computable sequence $(\mathcal{A}_{n})_{n\in\omega}$ such
that
\[
\mathcal{A}_{n}\cong\left\{
\begin{array}
[c]{ll}%
\mathcal{A}^{-} & \text{if}\ n\notin S_{1}\text{,}\\
\mathcal{A} & \text{if}\ n\in S_{1}-S_{2}\text{,}\\
\mathcal{A}^{+} & \text{if}\ n\notin S_{1}\cap S_{2}\text{.}%
\end{array}
\right.
\]
To accomplish this, let $D_{0}=D(\mathcal{A}^{-})$. At stage $s$, if $n\notin
S_{1,s}$, we let $D_{s}$ be the atomic diagram of $\mathcal{A^{-}}$. If $n\in
S_{1,s}-S_{2,s}$, we let $D_{s}$ be the atomic diagram of $\mathcal{A}$. If
$n\in S_{1,s}\cap S_{2,s}$, we let $D_{s}$ be the atomic diagram of
$\mathcal{A}^{+}$. There is some $s_{0}$ such that for all $s\geq s_{0}$,
$n\in S_{1}$ iff $n\in S_{1,s}$, and $n\in S_{2}$ iff $n\in S_{2,s}$. Let
$\mathcal{A}_{n}$ be the structure with diagram $D_{s}$ for $s\geq s_{0}$. It
is clear that $\mathcal{A}_{n}\cong\mathcal{A}$ iff $n\in S$.

\endproof

\section{Vector spaces}

\label{vssec}

The finite dimensional $\mathbb{Q}$-vector spaces over a fixed field are
completely determined by a finite set (a basis), so we might expect these to
behave much like finite structures. However, we have added complexity because
of the fact that for $1\leq m < n$, if $\mathcal{V}_{n}$ is a space of
dimension $n$, and $\mathcal{V}_{m}$ is an $m$-dimensional subspace, if
$\mathcal{V}_{n}\models\varphi(\overline{c},\overline{a})$, where $\varphi$ is
finitary quantifier-free, and $\overline{c}$ is in $\mathcal{V}_{m}$, then
there exists $\overline{a}^{\prime}$ such that $\mathcal{V}_{m}\models
\varphi(\overline{c},\overline{a}^{\prime})$. We work with vector spaces over
$\mathbb{Q}$, for concreteness, but any other infinite computable field would
give exactly the same results.

\begin{prop}
\label{prop3.1}

Let $K$ be the class of $\mathbb{Q}$-vector spaces, and let $\mathcal{A} $ be
a member of $K$.

\begin{enumerate}
\item If $\dim(\mathcal{A}) = 0$, then $I(\mathcal{A})$ is $m$-complete
$\Pi^{0}_{1}$ within $K$.

\item If $\dim(\mathcal{A}) = 1$, then $I(\mathcal{A})$ is $m$-complete
$\Pi^{0}_{2}$ within $K$.

\item If $\dim(\mathcal{A}) > 1$, then $I(\mathcal{A})$ is $m$-complete
$d$-$\Sigma^{0}_{2}$ within $K$.
\end{enumerate}
\end{prop}

\proof

For 1, first we note that $\mathcal{A}$ has a a finitary $\Pi_{1}$ Scott
sentence, within $K$, saying $(\forall x)\,x=0$. It follows that
$I(\mathcal{A})$ is $\Pi_{1}^{0}$ within $K$. Toward completeness, let $S$ be
a $\Pi_{1}^{0}$ set. We build a uniformly computable sequence of structures
$(\mathcal{A}_{n})_{n\in\omega}$ such that
\[
\dim(\mathcal{A}_{n})=\left\{
\begin{array}
[c]{ll}%
0 & \text{if}\ n\in S\text{,}\\
1 & \text{if}\ n\notin S\text{.}%
\end{array}
\right.
\]

Let $\mathcal{V}_{0}$ be a space of dimension $0$, and let $\mathcal{V}_{1}$
be a computable extension having dimension $1$. We have a computable sequence
$(S_{s})_{s\in\omega}$ of approximations for $S$ such that $n\in S$ iff for
all $s$, $n\in S_{s}$, and if $n\notin S_{s}$, then for all $t>s$, $n\notin
S_{t}$. If $n\in S_{s}$, we let $D_{s}=D(\mathcal{V}_{0})$. If $n\notin S_{s}%
$, then we let $D_{s}$ consist of the first $s$ sentences of $D(\mathcal{V}%
_{1})$. This completes the proof for 1.

\bigskip

Next, we turn to 2. First, we show that $\mathcal{A}$ has a computable
$\Pi_{2}$ Scott sentence. We have a computable $\Pi_{2}$ sentence
characterizing the class $K$. We take the conjunction of this with the
sentence saying
\[
(\exists x)\,\,x\neq0\ \&\ (\forall x)\,(\forall y)\,\bigvee\limits_{\lambda
\in\Lambda}\hspace{-0.15in}\bigvee\lambda(x,y)=0\text{,}%
\]
where $\Lambda$ is the set of all nontrivial linear combinations $q_{1}%
x+q_{2}y$, for $q_{i}\in\mathbb{Q}$. Now, $I(\mathcal{A})$ is $\Pi_{2}^{0}$.
We do not need to locate $\mathcal{A}$ within $K$, since the set of indices
for members of $K$ is $\Pi_{2}^{0}$.

For completeness, let $S$ be a $\Pi_{2}^{0}$ set. We build a uniformly
computable sequence $(\mathcal{A}_{n})_{n\in\omega}$ such that
\[
\dim(\mathcal{A}_{n})=\left\{
\begin{array}
[c]{ll}%
1 & \text{if}\ n\in S\text{,}\\
2 & \text{if}\ n\notin S\text{.}%
\end{array}
\right.
\]
We have computable approximations $S_{s}$ for $S$ such that $n\in S$ iff for
infinitely many $s$, $n\in S_{s}$.

Let $\mathcal{V}^{+}$ be a $2$-dimensional computable vector space and let
$\mathcal{V}$ be a\linebreak$1$-dimensional subspace. Let $B$ be an infinite
computable set of constants, for the universe of all $\mathcal{A}_{n}$. At
each stage $s$, we have a finite partial $1-1$ function $p_{s}$ from $B$ to a
target structure $\mathcal{V}$ (if $n\in S_{s}$) or $\mathcal{V}^{+}$ (if
$n\notin S_{s}$), and we have enumerated a finite set $D_{s}$ (a part of the
atomic diagram of $\mathcal{A}_{n}$) such that $p_{s}$ maps the constants
mentioned in $D_{s}$ into the target structure so as to make all of the
sentences true. We arrange that if $n\in S$, then $\cup\{p_{s}:$ $n\in
S_{s}\}$ maps $B$ onto $\mathcal{V}$. If $n\notin S$, then there is some
$s_{0}$ such that for all $s\geq s_{0}$, $n\notin S_{s}$. In this case,
$\cup_{s\geq s_{0}}p_{s}$ maps $B$ onto $\mathcal{V}^{+}$. For all $s$
(whether or not $n\in S_{s}$), $D_{s}$ decides the first $s$ atomic sentences
involving constants in $dom(p_{s})$.

We start with $p_{0}=\emptyset$, and $D_{0}=\emptyset$. Without loss of
generality, we may suppose that $n\in S_{0}$. We consider $p_{0}$ to be
mapping into $\mathcal{V}$. If there is no change in our guess about whether
$n\in S$ at stage $s+1$, then $p_{s+1}\supseteq p_{s}$, where the first $s+1$
constants from $B$ are in the domain, and the first $s+1$ elements of the
target structure are in the range. We must say what happens when we change our
guess at whether $n\in S$.

There are two cases. First, suppose $n\in S_{s+1}$ and $n\notin S_{s}$. In
this case, we take the greatest stage $t\leq s$ such that $n\in S_{t}$. We let
$p_{s+1}\supseteq p_{t}$ such that $p_{s+1}$ makes the sentences of $D_{s}$
true in $\mathcal{V}$, extending so that the first $s+1$ constants from $B$
are in the domain, and the first $s+1$ elements of $\mathcal{V}$ are in the
range. Next, assume $n\notin S_{s+1}$ and $n\in S_{s}$. In this case, we do
not look back at any earlier stage. We let $p_{s+1}\supseteq p_{s}$, extending
so that the first $s+1$ constants from $B$ are in the domain, and the first
$s+1$ elements of $\mathcal{V}^{+}$ are in the range. In either case, we let
$D_{s+1}\supseteq D_{s}$ so that for the first $s+1$ atomic sentences $\beta$
involving constants in the domain of $p_{s+1}$, $D_{s+1}$ includes $\pm\beta$,
whichever is made true by $p_{s+1}$. This completes the proof for 2.

\bigskip

Finally, we turn to 3. Suppose $\mathcal{A}$ has dimension $k$, where $k>1$.
Then $\mathcal{A}$ has a $d$-$\Sigma_{2}$ Scott sentence. We take the
conjunction of the axioms for $\mathbb{Q}$-vector spaces, and we add a
sentence saying that there are at least $k$ independent elements, and that
there are not at least $k+1$. Then $I(\mathcal{A})$ is $d$-$\Sigma_{2}^{0}$.
Toward completeness, let $S=S_{1}-S_{2}$, where $S_{1}$ and $S_{2}$ are both
$\Sigma_{2}^{0}$. Let $\mathcal{V}^{+}$ be a computable vector space of
dimension $k+1$, and let
\[
\mathcal{V}^{-}\subseteq\mathcal{V}\subseteq\mathcal{V}^{+}\text{,}%
\]
where $\mathcal{V}^{-}$ has dimension $k-1$, and $\mathcal{V}$ has dimension
$k$. We will produce a uniformly computable sequence of structures
$(\mathcal{A}_{n})_{n\in\omega}$ such that
\[
\dim(\mathcal{A}_{n})=\left\{
\begin{array}
[c]{ll}%
k-1 & \text{if}\ n\notin S_{1}\text{,}\\
k & \text{if}\ n\in S_{1}-S_{2}\text{,}\\
k+1 & \text{if}\ n\in S_{1}\cap S_{2}\text{.}%
\end{array}
\right.
\]

The construction is similar to that for 2. We have computable approximations
$S_{1,s}$ and $S_{2,s}$ for $S_{1}$ and $S_{2}$, respectively, such that
\[
n\in S_{i}\ \text{iff for all but finitely many }s\text{, }n\in S_{i,s}%
\text{.}%
\]
Let $B$ be an infinite computable set of constants, for the universe of all
$\mathcal{A}_{n}$. At each stage $s$, we have a finite partial $1-1$ function
$p_{s}$ from $B$ to $\mathcal{V}^{-}$ if $n\notin S_{1,s}$, to $\mathcal{V}$
if $n\in S_{1,s}-S_{2,s}$, and to $\mathcal{V}^{+}$ if $n\in S_{1,s}\cap
S_{2,s}$. We have $D_{s}$, a finite part of $D(\mathcal{A}_{n})$ such that
$p_{s}$ makes $D_{s}$ true in the target structure.

We arrange that if $n\notin S_{1}$, then the union of the $p_{s}$ for $n\notin
S_{1,s}$ maps $B$ onto $\mathcal{V}^{-}$. If $n\in S_{1}-S_{2}$, then at every
$s$ after some stage $s_{0}$, $n\in S_{1,s}$, while $n\in S_{2,s}$ for
infinitely many $s$. In this case, the union of $p_{s}$ for $s\geq s_{0}$ such
that $n\in S_{2,s}$ maps $B$ onto $\mathcal{V}$. If $n\in S_{1}\cap S_{2}$,
then for some stage $s_{1}$, for $s\geq s_{1}$, $n\in S_{1,s}\cap S_{2,s}$,
and the union of $p_{s}$ for $s\geq s_{1}$ maps $B$ onto $\mathcal{V}^{+}$.

We may suppose that $n\notin S_{1,0}$, so the target structure at stage $0$ is
$\mathcal{V}^{-}$. At stage $s+1$, if there is no change in the target
structure, then we extend $p_{s}$. We must say what to do when we change our
mind about the target structure. First, suppose the change is because of
$S_{1}$. If $n\in S_{1,s}$ and $n\notin S_{1,s+1}$, then the target structure
changes from $\mathcal{V}$ or $\mathcal{V}^{+}$ back to $\mathcal{V}^{-}$. We
take the greatest stage $t<s$ such that $n\notin S_{1,t}$. We let
$p_{s+1}\supseteq p_{t}$ such that $p_{s+1}$ makes $D_{s}$ true in
$\mathcal{V}^{-}$, extending so that the first $s+1$ constants from $B$ are in
the domain, and the first $s+1$ elements of $\mathcal{V}^{-}$ are in the
range. If $n\notin S_{1,s}$ and $n\in S_{1,s+1}-S_{2,s+1}$, then the target
structure changes from $\mathcal{V}^{-}$ to $\mathcal{V}$. We let
$p_{s+1}\supseteq p_{s}$, extending so that the first $s+1$ constants from $B$
are in the domain and the first $s+1$ elements of $\mathcal{V}$ are in the range.

Now, suppose the change is because of $S_{2}$. We suppose that $n\in S_{1,s}$
and $n\in S_{1,s+1}$. If $n\in S_{2,s}$ and $n\notin S_{2,s+1}$, then the
target structure changes from $\mathcal{V}^{+}$ back to $\mathcal{V}$. We take
the greatest stage $t\leq s$ such that the target structure is $\mathcal{V}$,
and we have had $n\in S_{1,s^{\prime}}$ for all $t<s^{\prime}<s$. If there is
no such $t$, then we take the greatest $t<s$ such that $n\notin S_{1,t}$. We
let $p_{s+1}\supseteq p_{t}$ such that $p_{s+1}$ makes $D_{s}$ true in
$\mathcal{V}$, extending to include the first $s+1$ elements of $B$ in the
domain and the first $s+1$ elements of $\mathcal{V}$ in the range. We let
$D_{s+1}\supseteq D_{s}$ so that for for the first $s+1$ atomic sentences
$\beta$ involving constants in the domain of $p_{s+1}$, $D_{s+1}$ includes
$\pm\beta$, whichever is made true by $p_{s+1}$. This completes the proof of 3.

\endproof

The final case, that of an infinite dimensional space, is already known
\cite{C1}. We include it here for completeness.

\begin{prop}
Let $K$ be the class of computable vector spaces over $\mathbb{Q}$, and let
$\mathcal{A}$ be a member of $K$ of infinite dimension. Then $I(\mathcal{A})$
is $m$-complete $\Pi_{3}^{0}$ within $K$.
\end{prop}

\proof

We have a computable $\Pi_{3}$ Scott sentence for $\mathcal{A}$, obtained by
taking the conjunction of the axioms for $\mathbb{Q}$-vector spaces and the
conjunction over all $k\in\omega$ of computable $\Sigma_{2}$ sentences saying
that the dimension is at least $k$. Therefore, $I(\mathcal{A})$ is $\Pi
_{3}^{0}$.

For completeness, let $Cof$ denote the set of indices for cofinite c.e.\ sets.
It is well known that the complement of $Cof$ is $m$-complete $\Pi_{3}^{0}$
(see \cite{soare}). We build a uniformly computable sequence of vector spaces
$(\mathcal{A}_{n})_{n\in\omega}$ such that $\mathcal{A}_{n}$ has infinite
dimension iff $n\notin Cof$.

Let $\mathcal{V}$ be an infinite dimensional vector space with basis
$\{v_{i}:i\geq-1\}$. Let $B$ be an infinite computable set of constants, for
the universe of all $\mathcal{A}_{n}$. For each set $S\subseteq\omega$, let
$\mathcal{V}_{S}$ be the linear span of $\{v_{-1}\}\cup\{v_{i}:i\in S\}$. Our
goal is to make $\mathcal{A}_{n}\cong\mathcal{V}_{S}$, where $S=\omega-W_{n}$.
At stage $s$, we have a finite approximation of $S$, as follows. Let
$S_{0}=\emptyset$. If $W_{n+1,s+1}$ includes some $x\in S_{s}$, we let
$S_{s+1}$ be the result of removing from $S_{s}$ all $y\geq x$. If
$W_{n+1,s+1}$ contains no elements of $S_{s}$, we let $S_{s+1}$ be the result
of adding to $S_{s}$ the first element of $\omega$ not in $W_{n+1,s+1}$. Note
that for each $k$, there exists $s$ such that for all $t\geq s$, $S\cap
k=S_{t}\cap k$. Moreover, for every $s$, there is some $t\geq s$ such that
$S_{t}\subseteq\omega-W_{n}$. For such $t$, for all $t^{\prime}\geq t$,
$S_{t^{\prime}}\supseteq S_{t}$.

For the construction, at stage $s$ we have a finite partial $1-1$ function
$p_{s}$ from $B$ into $\mathcal{V}_{S_{s}}$. We include the first $s$ elements
of $B$ in the domain, and the first $s$ elements of $\omega$ that are in
$\mathcal{V}_{S_{s}}$ in the range. We also have $D_{s}$ deciding the first
$s$ atomic sentences with constants in $dom(p_{s})$, such that $p_{s}$ makes
the sentences true in $\mathcal{V}_{S_{s}}$. We start with $p_{0}=\emptyset$,
$D_{0}=\emptyset$, and we think of $p_{0}$ as mapping into $\mathcal{V}%
_{S_{0}}$.

At stage $s+1$, we define $p_{s+1}$ as follows. First, suppose $S_{s+1}$ is
the result of adding an element to $S_{s}$. Then $p_{s+1}\supseteq p_{s}$
including the first $s+1$ elements of $B$ in the domain, and the first $s+1$
elements of $\omega$ that are in $\mathcal{V}_{s+1}$ in the range. Now,
suppose $S_{s+1}$ is the result of removing one or more elements from $S_{s}$,
so that $S_{s+1}=S_{t}$ for some greatest $t<s$. We take $p_{s+1}\supseteq
p_{t}$ such that $p_{s+1}$ makes $D_{s}$ true in $\mathcal{V}_{s+1}%
=\mathcal{V}_{t}$. We let $D_{s+1}$ extend $D_{s}$ so as to decide the first
$s+1$ atomic sentences involving constants in $dom(p_{s+1})$. This completes
the construction.

There is an infinite sequence of stages $t$ such that for all $s>t$,
$p_{s}\supseteq p_{t}$. Let $f$ be the union of the functions $p_{t}$ for
these $t$. We can see that $f$ is a $1-1$ mapping of $B$ onto $\mathcal{V}%
_{S}$. Moreover, if $\mathcal{A}_{n}\cong_{f}\mathcal{V}_{S}$, then $\cup
_{s}D_{s}$ is the atomic diagram of $\mathcal{A}_{n}$.

\endproof

We have completely characterized the $m$-degrees of index sets of computable
vector spaces over $\mathbb{Q}$, within the class of all such vector spaces.

\section{Archimedean ordered fields}

\label{aofsec}

Archimedean ordered fields are isomorphic to subfields of the reals. They are
determined by the Dedekind cuts that are filled. It follows that for any
computable Archimedean ordered field, the index set is $\Pi^{0}_{3}$. In some
cases, the index set is simpler. We consider Archimedean real closed ordered
fields. We show that the index set is $m$-complete $\Pi^{0}_{2}$ if there are
no transcendentals, $m$-complete $d$-$\Sigma^{0}_{2}$ if the transcendence
degree is finite but not $0$, and $m$-complete $\Pi^{0}_{3}$ if the
transcendence degree is infinite.

\begin{prop}
\label{prop4.1}

If $\mathcal{A}$ is a computable Archimedean ordered field, then
$I(\mathcal{A})$ is $\Pi^{0}_{3}$.
\end{prop}

\proof

It is enough to show that $\mathcal{A}$ has a computable $\Pi_{3}$ Scott
sentence. We have a computable $\Pi_{2}$ sentence $\sigma_{0}$ characterizing
the Archimedean ordered fields. For each $a\in\mathcal{A}$, we have a
computable $\Pi_{1}$ formula $c_{a}(x)$ saying that $x$ is in the cut
corresponding to $a$---we take the conjunction of a c.e.\ set of formulas
saying $q<x<r$, for rationals $q,r$ such that $\mathcal{A}\models q<a<r$. Let
$\sigma_{1}$ be $\bigwedge\hspace{-0.18in}\bigwedge\,_{a}\,(\exists
x)\,c_{a}(x)$, and let $\sigma_{2}$ be $(\forall x)\,\bigvee\hspace
{-0.18in}\bigvee\,_{a}c_{a}(x)$. The conjunction of $\sigma_{0}$, $\sigma_{1}%
$, and $\sigma_{2}$ is a Scott sentence, which we may take to be computable
$\Pi_{3}$.

\endproof

We turn to real closed fields. Here is the main result of the section.

\begin{thm}
\label{thm4.2}

Let $\mathcal{A}$ be a computable Archimedean real closed ordered field.

\begin{enumerate}
\item If the transcendence degree is $0$ $($i.e., $\mathcal{A}$ is isomorphic
to the ordered field of algebraic reals$)$, then $I(\mathcal{A})$ is
$m$-complete $\Pi_{2}^{0}$ $($within the class of Archimedean real closed
ordered fields$)$.

\item If the transcendence degree is finite but not $0$, then $I(\mathcal{A})$
is $m$-complete $d$-$\Sigma_{2}^{0}$ $($within the class of Archimedean real
closed ordered fields$)$.

\item If the transcendence degree is infinite, then $I(\mathcal{A})$ is
$m$-complete $\Pi_{3}^{0}$ $($within the class of Archimedean real closed
ordered fields$)$.
\end{enumerate}
\end{thm}

It is convenient to have the following definition.

\begin{dfn}
Let $\mathcal{C}$ be a substructure of $\mathcal{B}$. We write $\mathcal{C}%
\preceq_{1}\mathcal{B}$ if for all finitary quantifier-free formulas $\varphi
$, if $\mathcal{B}\models\varphi(\overline{c},\overline{b})$, where
$\overline{c}$ is in $\mathcal{C}$, then there exists $\overline{b}^{\prime}$
such that $\mathcal{C}\models\varphi(\overline{c},\overline{b}^{\prime})$
$($in other words, satisfaction of finitary existential formulas by tuples in
$\mathcal{C}$ is the same in $\mathcal{C}$ and $\mathcal{B})$.
\end{dfn}

\noindent\textbf{Note}: We are using the notation $\preceq_{1}$ to avoid
confusion with $\leq_{1}$, which is often used for a weaker
relation---$\mathcal{C}\leq_{1}\mathcal{B}$ if the existential sentences true
in $\mathcal{B}$ are all true in $\mathcal{C}$.

\bigskip

We shall use the following lemma for Part 1 of Theorem \ref{thm4.2}.

\begin{lem}
\label{QaLemma}

Let $\mathcal{A}$ be the ordered field of real algebraic numbers, and let
$\mathcal{B}$ be the real closure of $e$ $($the familiar real
transcendental$)$. Then $\mathcal{A}\preceq_{1}\mathcal{B}$.
\end{lem}

\proof

Let $\varphi$ be finitary quantifier-free. Suppose $\mathcal{B}\models
\varphi(\overline{c},e,\overline{a})$, where $\overline{c}$ is in
$\mathcal{A}$. There is an open interval $I$ (in $\mathbb{R}$) containing $e$
such that for all $e^{\prime}\in I$, there exists $\overline{a}^{\prime}$ such
that $\mathbb{R}\models(\exists\overline{u})\,\varphi(\overline{c},e^{\prime
},\overline{a}^{\prime})$, where $\overline{a}^{\prime}$ is algebraic over
$\overline{c},e^{\prime}$ in the same way that $\overline{a}$ is algebraic
over $\overline{c},e$. Taking $e^{\prime}$ to be rational in $I$, we have
$\overline{a}^{\prime}$ in $\mathcal{A}$ such that $\mathcal{A}\models
\varphi(\overline{c},e^{\prime},\overline{a}^{\prime})$.

\endproof

\noindent\emph{Proof} of Part 1 of Theorem \ref{thm4.2}: Let $\mathcal{A}$ be
the ordered field of real algebraic numbers. We must show that $I(\mathcal{A}%
)$ is $m$-complete $\Pi_{2}^{0}$ (within the class of real closed Archimedean
ordered fields). To show that $I(\mathcal{A})$ is $\Pi_{2}^{0}$, we note that
there is a computable $\Pi_{2}$ Scott sentence. We take the conjunction of a
sentence characterizing the real closed ordered fields, and a sentence saying
that each element is a root of some polynomial.

We now show $m$-completeness. Let $S$ be a $\Pi_{2}^{0}$ set. Let
$\mathcal{B}$ be the real closure of $e$, as in Lemma \ref{QaLemma}. We
produce a uniformly computable sequence $(\mathcal{A}_{n})_{n\in\omega}$ such
that
\[
\mathcal{A}_{n}\cong\left\{
\begin{array}
[c]{ll}%
\mathcal{A} & \text{if }n\in S\text{,}\\
\mathcal{B} & \text{otherwise.}%
\end{array}
\right.
\]
Let $\mathcal{M}^{+}$ be a computable copy of $\mathcal{B}$, and let
$\mathcal{M}$ be the substructure isomorphic to $\mathcal{A}$. What is
important is that $\mathcal{M}^{+}$ is computable, $\mathcal{M}$ is a proper
c.e.\ substructure, and $\mathcal{M}\preceq_{1}\mathcal{M}^{+}$. Let $B$ be an
infinite computable set of constants, for the universe of all $\mathcal{A}%
_{n}$.

We have a computable approximation $(S_{s})_{s\in\omega}$ for $S$ such that
\[
n\in S\text{ iff }n\in S_{s}\ \text{for infinitely many }s.
\]
We construct $\mathcal{A}_{n}$ in stages. At each stage, we determine a finite
partial $1-1$ function $f_{s}$ from $B$ to the target structure, which is
$\mathcal{M}$ if $n\in S_{s}$, and $\mathcal{M}^{+}$ otherwise. The domain of
$f_{s}$ includes the first $s$ constants from $B$, and the range includes the
first $s$ elements of the target structure. Also, at stage $s$, we enumerate a
finite part of $D_{s}$ of the diagram of $\mathcal{A}_{n}$, such that $f_{s}$
makes $D_{s}$ true in the target structure. For the first $s$ atomic sentences
$\varphi$, involving only the first $s$ constants, we put $\pm\varphi$ in
$D_{s}$. We arrange that if $n\in S$, then $f=\cup\{f_{s}:$ $n\in S_{s}\}$ is
an isomorphism from $\mathcal{A}_{n}$ onto $\mathcal{M}_{1}$, and if $n\notin
S$, and $s_{0}$ is least such that for all $s\geq s_{0}$, $n\notin S_{s}$,
then $f=\cup_{s\geq s_{0}}f_{s}$ is an isomorphism from $\mathcal{A}_{n}$ onto
$\mathcal{M}^{+}$.

At stage $0$, we let $f_{0}=\emptyset$, and $D_{0}=\emptyset$. At stage $s+1$,
if $n\in S_{s+1}$ and $n\in S_{s}$, or if $n\notin S_{s+1}$ and $n\notin
S_{s}$, then we let $f_{s+1}\supseteq f_{s}$, adding to the domain and range.
We extend $D_{s}$ to $D_{s+1}$ so that $f_{s+1}$ makes the sentences of
$D_{s+1}$ true in the target structure. Suppose $n\in S_{s+1}$ and $n\notin
S_{s}$. Let $t<s$ be greatest such that $n\in S_{t}$ or $t=0$. We may suppose
that $f_{s}$ extends $f_{t}$. It follows from Lemma~\ref{QaLemma} that there
is an extension of $f_{t}$ which makes $D_{s}$ true in $\mathcal{M}$. We let
$f_{s+1}$ be such an extension that also includes the required elements in the
domain and range. We extend $D_{s}$ to $D_{s+1}$ so that $f_{s+1}$ makes the
sentences of $D_{s+1}$ true. Finally, suppose $n\in S_{s}$ and $n\notin
S_{s+1}$. We let $f_{s+1}\supseteq f_{s}$, changing to the larger target
structure. We add the required elements to the domain and range. We extend
$D_{s}$ to $D_{s+1}$ so that $f_{s+1}$ makes the sentences of $D_{s+1}$ true
in the target structure. This completes the proof of Part 1.

\bigskip\noindent\emph{Proof} of Part 2 of Theorem \ref{thm4.2}: Let
$\mathcal{A}$ be the real closure of $\{e_{1},\ldots,e_{k}\}$, where the
$e_{i}$ are algebraically independent computable reals, $k\geq1$. We must show
that $I(\mathcal{A})$ is $m$-complete $d$-$\Sigma_{2}^{0}$ (within the class
of Archimedean real closed ordered fields). To show that $I(\mathcal{A})$ is
$d$-$\Sigma_{2}$, note that there is a computable $\Sigma_{2}$ sentence
$\varphi$ saying that there are elements $x_{1}^{\mathcal{A}},\dots
,x_{k}^{\mathcal{A}}$ filling the cuts of $e_{1},\dots,e_{k}$. There is
another computable $\Sigma_{2}$ sentence $\psi$ saying that there are elements
$x_{1},\ldots,x_{k}$ filling these cuts, and another element $y$ not equal to
any rational function of the $x_{i}$. Then the conjunction of $\varphi
\ \&\ \lnot{\psi}$ with the computable $\Pi_{2}$ sentence characterizing the
Archimedean real closed ordered fields is a Scott sentence for $\mathcal{A}$.
It follows that $I(\mathcal{A})$ is $d$-$\Sigma_{2}^{0}$.

Toward completeness, let $S=S_{1}-S_{2}$, where $S_{i}$ is $\Sigma_{2}^{0}$.
Let $e_{k+1}$ be a computable real number algebraically independent of
$e_{1},\dots,e_{k}$. Let $\mathcal{M}^{+}$ be a computable structure
isomorphic to the real closure of $\{e_{1},\ldots,e_{k},e_{k+1}\}$, let
$\mathcal{M}$ be the subfield isomorphic to $\mathcal{A}$, and let
$\mathcal{M}^{-}$ be the subfield generated by the elements corresponding to
$e_{i}$ for $i < k$. (If $k = 1$, then $\mathcal{M}^{-}$ is isomorphic to the
ordered field of real algebraic numbers.) We have the following lemma, proved
in the same way as Lemma \ref{QaLemma}.

\begin{lem}
$\mathcal{M}^{-}\preceq_{1}\mathcal{M}\preceq_{1}\mathcal{M}^{+}$
\end{lem}

What is important is that $\mathcal{M}^{+}$ is a computable structure,
$\mathcal{M}^{-}$ and $\mathcal{M}$ are c.e.\ substructures, and the three
structures form a properly increasing chain, with $\mathcal{M}^{-}\preceq
_{1}\mathcal{M}\preceq_{1}\mathcal{M}^{+}$. We will produce a uniformly
computable sequence $(\mathcal{A}_{n})_{n\in\omega}$ such that
\[
\mathcal{A}_{n}\cong\left\{
\begin{array}
[c]{ll}%
\mathcal{M}^{-} & \text{if }n\notin S_{1}\text{,}\\
\mathcal{M} & \text{if }n\in S_{1}-S_{2}\text{,}\\
\mathcal{M}^{+} & \text{if }n\in S_{1}\cap S_{2}\text{.}%
\end{array}
\right.
\]

Let $B$ be an infinite computable set of constants, for the universe of all
$\mathcal{A}_{n}$. We have computable approximations $(S_{i,s})_{s\in\omega}$
for $S_{i}$ such that
\[
n\in S_{i}\ \text{iff }(\exists s_{0})\,(\forall s\geq s_{0})\,[n\in
S_{i,s}]\text{.}%
\]
At each stage $s$, we will have a finite partial $1-1$ function $f_{s}$ from
$B$ to the appropriate target structure, which is $\mathcal{M}^{-}$ if
$n\notin S_{1,s}$, $\mathcal{M}$ if $n\in S_{1,s}-S_{2,s}$, and $\mathcal{M}%
^{+}$ if $n\in S_{1,s}\cap S_{2,s}$. The domain of $f_{s}$ will include the
first $s$ constants from $B$, and the range will include the first $s$
elements of the target structure. We will also have a finite set $D_{s}$ of
atomic sentences and negations of atomic sentences such that $f_{s}$ makes
$D_{s}$ true in the target structure. For the first $s$ atomic sentences
$\varphi$ involving only the first $s$ constants, $D_{s}$ will include
$\pm\varphi$. We let $\mathcal{A}_{n}$ be the structure with atomic diagram
$\cup_{s}D_{s}$.

We will arrange that if $n\notin S_{1}$, then the union of $f_{s}$ for
$n\notin S_{1,s}$ will be an isomorphism from $\mathcal{A}_{n}$ onto
$\mathcal{M}^{-}$. If $n\in S_{1} - S_{2}$, and $s_{0}$ is first such that for
all $s\geq s_{0}$, $n\in S_{1,s}$, then the union of $f_{s}$, for $s\geq
s_{0}$ such that $n\notin S_{2,s} $, will be an isomorphism from
$\mathcal{A}_{n}$ onto $\mathcal{M}$. If $n\in S_{1}\cap S_{2}$, and $s_{1}$
is first such that for all $s\geq s_{1}$, $n\in S_{1,s} $ and $n\in S_{2,s}$,
then the union of $f_{s}$ for $s\geq s_{1}$ will be an isomorphism from
$\mathcal{A}_{n}$ onto $\mathcal{M}^{+}$.

We start with $f_{0}=\emptyset$ and $D_{0}=\emptyset$. Suppose we have $f_{s}$
and $D_{s}$. If the target structure at stage $s+1$ is the same as at stage
$s$, then we extend $f_{s}$ and $D_{s}$ in the obvious way. Suppose $n\notin
S_{1,s+1}$, where $n\in S_{1,s}$. Let $t<s$ be greatest such that $n\notin
S_{1,t}$ or $t=0$. Say $f_{t}$ maps $\overline{d}$ to $\overline{c}$ in
$\mathcal{M}^{-}$. We let $f_{s+1}\supseteq f_{t}$, where $f_{s+1}$ makes the
sentences of $D_{s}$ true in $\mathcal{M}^{-}$. Now, suppose $n\in
S_{1,s+1}-S_{2,s+1}$, where $n\in S_{1,s}\cap S_{2,s}$. Take greatest $t<s$
such that $n\in S_{1,t}-S_{2,t}$ and for all $t^{\prime}$ with $t<t^{\prime
}<s$, $n\in S_{1,t^{\prime}}$. Then we let $f_{s+1}\supseteq f_{t}$, where
$f_{s+1}$ makes the sentences of $D_{s}$ true in $\mathcal{M}$.

Suppose that $n\in S_{1,s+1}-S_{2,s+1}$, where $n\notin S_{1,s}$. Then we take
$f_{s+1}\supseteq f_{s}$, where $f_{s+1}$ makes $D_{s}$ true in the target
structure $\mathcal{M}$. Similarly, if we have $n\in S_{1,s+1}\cap S_{2,s+1}$,
where either $n\notin S_{1,s}$ or $n\in S_{1,s}-S_{2,s}$, we let
$f_{s+1}\supseteq f_{s}$, where $f_{s+1}$ makes $D_{s}$ true in the target
structure $\mathcal{M}^{+}$. We extend $f_{s}$ to include the required
elements in the domain and range, and we let $D_{s+1}\supseteq D_{s} $ so that
$f_{s+1}$ makes $D_{s+1}$ true. This completes the proof of Part 2.

\bigskip\noindent\emph{Proof} of Part 3 of Theorem \ref{thm4.2}: Let
$\mathcal{A}$ be a computable real closed Archimedean ordered field of
infinite transcendence degree. We must show that $I(\mathcal{A})$ is
$m$-complete $\Pi_{3}^{0}$ (within the class of real closed Archimedean
ordered fields). It follows from Proposition \ref{prop4.1} that $I(\mathcal{A}%
)$ is $\Pi_{3}^{0}$. For completeness, we need the following lemma.

\begin{lem}
\label{lem4.6}

There is a uniformly c.e.\ sequence $(\mathcal{M}_{k})_{k\in\omega}$ of real
closed subfields of $\mathcal{A}$, such that for all $k$, $\mathcal{M}_{k+1}$
properly extends $\mathcal{M}_{k}$, $\mathcal{M}_{k}\preceq_{1}\mathcal{M}%
_{k+1}$, and $\cup_{k}\mathcal{M}_{k} = \mathcal{A}$.
\end{lem}

\proof

Let $(a_{n})_{n\in\omega}$ be a computable list of all elements of
$\mathcal{A}$. Let $n(k)$ be the least $n$ such that $a_{n}$ is not algebraic
over $\{a_{n(j)}:$ $j<k\}$. Let $\mathcal{M}_{k}$ be the real closure of
$\{a_{n(j)}:$ $j<k\}$ in $\mathcal{A}$, which is the same as the real closure
of $\{a_{n}:$ $n<n(k)\}$. It is clear that $\mathcal{M}_{k+1}$ properly
extends $\mathcal{M}_{k}$, and $\mathcal{A}=\cup_{k}\mathcal{M}_{k}$. We have
a computable function $f(k,s)$, which, for each $k$, is nondecreasing, with
limit $n(k)$. Using this, we can see that $\mathcal{M}_{k}$ is c.e., uniformly
in $k$.

We must show that $\mathcal{M}_{k}\preceq_{1}\mathcal{M}_{k+1}$. Suppose that
$\mathcal{M}_{k+1}\models\varphi(\overline{c},\overline{b})$, where $\varphi$
is quantifier-free, and $\overline{c}$ is in $\mathcal{M}_{k}$. We may suppose
that $\overline{c}$ includes $a_{m}$ for all $m<n(k)$, and $\overline{b}$ has
the form $a_{n_{k}},\overline{d}$. Say $\psi(\overline{c},a_{n(k)}%
,\overline{d})$ states how $\overline{d}$ are expressed as roots of
polynomials over functions of $\overline{c},a_{n(k)}$. There is an interval
around $a_{n(k)}$, with rational endpoints, such that for all $x\in I$, and
all $\overline{u}$ satisfying $\psi(\overline{c},x,\overline{u})$,
$\varphi(\overline{c},x,\overline{u})$ holds in $\mathbb{R}$. Taking $x$
rational in this interval, we have $\overline{u}$ satisfying $\psi
(\overline{c},x,\overline{u})$ in $\mathcal{M}_{k}$, so $\mathcal{M}%
_{k}~\models~\varphi(\overline{c},x,\overline{u})$.

\endproof

Recall that $Cof=\{n:$ $\omega-W_{n}\ $is finite$\}$. This set is $m$-complete
$\Sigma_{3}^{0}$, so the complement is $m$-complete $\Pi_{3}^{0}$. We have a
$\Delta_{2}^{0}$ function $\nu(n,r)$, nondecreasing in $r$, for each $n$, such
that if $n\in Cof$, then $\lim_{r}\nu(n,r)$ has value equal to the cardinality
of $\omega-W_{n}$, and if $n\notin Cof$, then $\lim_{r}\nu(n,r)=\infty$. We
have a computable approximation to $\nu$, and we define a computable function
$g(n,s)$, such that $g(n,0)$ is our stage $0$ guess at $\nu(n,0)$. Supposing
that $g(n,s)$ is our stage $s$ guess at $\nu(n,r)$, if our stage $s+1$ guess
at $\nu$ agrees with the stage $s$ guess at $\nu(n,j)$ for all $j\leq r$, then
$g(n,s+1)$ is the stage $s+1$ guess at $\nu(n,r+1)$. If there is a
disagreement, say $j$ is the least such that our stage $s$ and stage $s+1$
guesses at $\nu(n,j)$ disagree. Let $g(n,s+1)$ be the stage $s+1$ guess at
$\nu(n,j)$. If $n\notin Cof$, then for all $k$, there exists $s_{0}$ such
that
\[
(\forall s\geq s_{0})\,[g(n,s)\geq k]\text{.}%
\]
If $n\in Cof$, and $k$ is the cardinality of $\omega-W_{n}$, then there are
infinitely many $s$ such that $g(n,s)=k$, and $k$ is the least such number. In
other words, $\liminf_{s}g(n,s)$ is the cardinality of $\omega-W_{n}$.

We will build a uniformly computable sequence $(\mathcal{A}_{n})_{n\in\omega}$
of Archimedean ordered fields such that if $n\notin Cof$, that is,
$\omega-W_{n}$ is infinite, then $\mathcal{A}_{n}\cong\mathcal{A}$. If $n\in
Cof$, that is, $\omega-W_{n}$ has size $k$ for some $k\in\omega$, then
$\mathcal{A}_{n}\cong\mathcal{M}_{k}$. Let $B$ be an infinite computable set
of constants, for the universe of all $\mathcal{A}_{n}$. At stage $s$, we have
a finite partial $1-1$ function $f_{s}$ from $B$ to the target structure
$\mathcal{M}_{k}$, where $k=g(n,s)$. The domain of $f_{s}$ will include the
first $s$ constants from $B$, and the range will include the first $s$
elements of the target structure. We also have $D_{s}$, a finite set of atomic
sentences and negations of atomic sentences which $f_{s}$ makes true in the
target structure. For each $\varphi$ among the first $s$ atomic sentences
using only the first $s$ constants, $D_{s}$ will include $\pm\varphi$.

We let $\mathcal{A}_{n}$ be the structure with
\[
D(\mathcal{A}_{n})=\cup_{s}D_{s}\text{.}%
\]
Let $T$ be the set of $s$ such that for all $t\geq s$, $g(n,t)\geq g(n,s)$. We
shall arrange that $f=\cup_{s\in T}f_{s}$ is an isomorphism from
$\mathcal{A}_{n}$ onto the desired structure. With this goal in mind, we
maintain the following condition.

\bigskip\noindent\textbf{Condition maintained}: Suppose $t<s$, where
$g(n,t)\leq g(n,s)$, and for all $t^{\prime}$ with $t<t^{\prime}<s$, we have
$g(n,t^{\prime})\geq g(n,t)$. Then $f_{s}\supseteq f_{t}$.

\bigskip

Let $f_{0}=\emptyset$, and let $D_{0}=\emptyset$. Given $f_{s}$ and $D_{s}$,
we must determine $f_{s+1}$ and $D_{s+1}$. The target structure is
$\mathcal{M}_{k}$, where $k=g(n,s+1)$. First, suppose $g(n,s+1)\geq g(n,s)$.
Then we let $f_{s+1}\supseteq f_{s}$. We extend to include the required
elements in the domain and range. We let $D_{s+1}\supseteq D_{s}$, where
$f_{s+1}$ makes the sentences true in the target structure. Now, suppose
$g(n,s+1)<g(n,s)$. Let $t<s$ be greatest such that $g(n,t)=g(n,s+1)$ and there
is no $t^{\prime}$ with $t<t^{\prime}<s$ such that $g(n,t^{\prime})<g(n,t)$,
or if there is no $t<s$ such that $g(n,t)=g(n,s+1)$, take the greatest $t<s$
such that $g(n,t)\leq g(n,s+1)$. We may assume that $f_{s}$ extends $f_{t}$.
We let $f_{s+1}$ extend $f_{t}$ such that $f_{s+1}$ makes $D_{s}$ true in the
target structure. This completes the proof of Part 3. So, we have finished the
proof of Theorem \ref{thm4.2}.

\endproof

\subsection{Further remarks}

We can say a little about Archimedean ordered fields that are not real closed,
but the results are fragmentary. We mention two results, without proof. First,
suppose $\mathcal{A}$ is a computable Archimedean ordered field which is a
purely transcendental extension of $\mathbb{Q}$. If the transcendence degree
is $0$ (i.e., $\mathcal{A}$ is isomorphic to the ordered field of rationals),
then $I(\mathcal{A})$ is $m$-complete $\Pi_{2}^{0}$ (within the class of
Archimedean ordered fields). If the transcendence degree is finite but not
$0$, then $I(\mathcal{A})$ is $m$-complete $d$-$\Sigma_{2}^{0}$ (within the
class of Archimedean ordered fields). We do not have a result for infinite
transcendence degree. Next, suppose $\mathcal{A}$ is a computable Archimedean
ordered field that is a finite algebraic extension of $\mathbb{Q}$. Then,
again, $I(\mathcal{A})$ is $m$-complete $\Pi_{2}^{0}$. If we extend further,
adding finitely many algebraically independent computable reals, then the
index set is $m$-complete $d$-$\Sigma_{2}^{0}$.

\section{Abelian $p$-groups}

\label{apgsec}

\subsection{Preliminaries on Abelian $p$-groups}

For a prime $p$, a $p$\emph{-group} is a group in which every element has
order $p^{n}$ for some $n$. We will consider only countable reduced Abelian
$p$-groups. These groups are of particular interest because of their
classification up to isomorphism by Ulm. For analysis of classical Ulm's
Theorem and a more detailed discussion of this class of groups, consult
Kaplansky's book \cite{kaplansky}. Generally, our notation here will be
similar to Kaplansky's.

Let $G$ be a countable Abelian $p$-group. We define a sequence of subgroups
$G_{\alpha}$, letting $G_{0}=G$, $G_{\alpha+1}=pG_{\alpha}$, and for limit
$\alpha$, $G_{\alpha}=\cap_{\beta<\alpha}G_{\beta}$. There is a countable
ordinal $\alpha$ such that $G_{\alpha}=G_{\alpha+1}$. The least such $\alpha$
is the \emph{length} of $G$, denoted by $\lambda(G)$. The group is
\emph{reduced} if $G_{\lambda(G)}=\{0\}$. An element $x\not =0$ has
\emph{height} $\beta$ if $x\in G_{\beta}-G_{\beta+1}$. Let $P(G)$ be the set
of element of $G$ of order $p$. Let $P_{\alpha}=G_{\alpha}\cap P(G)$. For each
$\beta<\lambda(G)$, $P_{\beta}/P_{\beta+1}$ is a vector space over
$\mathbb{Z}_{p}$ of dimension $\leq\aleph_{0}$, and this dimension is denoted
by $u_{\beta}(G)$. The \emph{Ulm sequence} for $G$ is the sequence $(u_{\beta
}(G))_{\beta<\lambda(G)}$.

For any computable ordinal $\alpha$, it is fairly straightforward to write a
computable infinitary sentence stating that $G$ is a reduced Abelian $p$-group
of length at most $\alpha$, and describing its Ulm invariants. In particular,
Barker \cite{B} established the following results.

\begin{lem}
\label{sentapg}

Let $G$ be a computable Abelian $p$-group.

\begin{enumerate}
\item $G_{\omega\cdot\alpha}$ is $\Pi^{0}_{2 \alpha}$.

\item $G_{\omega\cdot\alpha+ m}$ is $\Sigma^{0}_{2 \alpha+ 1}$.

\item $P_{\omega\cdot\alpha}$ is $\Pi^{0}_{2 \alpha}$.

\item $P_{\omega\cdot\alpha+ m}$ is $\Sigma^{0}_{2 \alpha+ 1}$.
\end{enumerate}
\end{lem}

\proof

It is easy to see that 3 and 4 follow from 1 and 2, respectively. Toward 1 and
2, note the following:%

\begin{align*}
x\in G_{m}  &  \iff\exists y[p^{m}y=x]\text{;}\\
x\in G_{\omega}  &  \iff\bigwedge\limits_{m\in\omega}\hspace{-0.15in}%
\bigwedge\exists y[p^{m}y=x]\text{;}\\
x\in G_{\omega\cdot\alpha+m}  &  \iff\exists y[p^{m}y=x\ \&\ G_{\omega
\cdot\alpha}(y)]\text{;}\\
x\in G_{\omega\cdot\alpha+\omega}  &  \iff\bigwedge\limits_{m\in\omega}%
\hspace{-0.15in}\bigwedge\exists y[p^{m}y=x\ \&\ G_{\omega\cdot\alpha
}(y)]\text{;}\\
x\in G_{\omega\cdot\alpha}  &  \iff\bigwedge\limits_{\gamma<\alpha}%
\hspace{-0.15in}\bigwedge G_{\omega\cdot\gamma}(x)\text{ for limit\ }%
\alpha\text{.}%
\end{align*}

\endproof

Using Lemma \ref{sentapg}, it is easy to write, for any computable ordinal
$\beta$, a computable $\Pi_{2\beta+1}^{0}$ sentence whose models are exactly
the reduced Abelian $p$-groups of length $\omega\beta$.

Khisamiev \cite{kh4}, \cite{kh3} gave a very good characterization of the
reduced Abelian $p$-groups of length less than $\omega^{2}$, which have
computable isomorphic copies. For groups of finite length, it is easy to
produce computable copies. Khisamiev gave a characterization for length
$\omega$, and proved an inductive lemma that allowed him to build up to all
lengths less than $\omega^{2}$. Here is the result for length $\omega$.

\begin{prop}
[Khisamiev]\label{Kh1}

Let $\mathcal{A}$ be a reduced Abelian $p$-group of length~$\omega$. Then
$\mathcal{A}$ has a computable copy iff

\begin{enumerate}
\item the relation $R_{\mathcal{A}}=\{(n,k):$ $u_{n}(\mathcal{A})\geq k\}$ is
$\Sigma_{2}^{0}$, and

\item there is a computable function $f_{\mathcal{A}}$ such that for each $n$,
$f_{\mathcal{A}}(n,s)$ is nondecreasing and with limit $n^{\ast}\geq n$ such
that $u_{n^{\ast}}(\mathcal{A})\not =0$.
\end{enumerate}

\noindent Moreover, we can effectively determine a computable index for a copy
of $\mathcal{A}$ from a $\Sigma_{2}^{0}$ index for $R_{\mathcal{A}}$, and a
computable index for a function~$f_{\mathcal{A}}$.
\end{prop}

Here is the inductive lemma.

\begin{lem}
[Khisamiev]\label{Kh2}

Let $\mathcal{A}$ be a reduced Abelian $p$-group. Suppose $\mathcal{A}%
_{\omega}$ is $\Delta_{3}^{0}$, the relation $R_{\mathcal{A}}$ is $\Sigma
_{2}^{0}(X)$, and there is a function $f_{\mathcal{A}}$ such that for all $n$,
$f_{\mathcal{A}}(n,s)$ is nondecreasing and with limit $r^{\ast}\geq r$ such
that $u_{r^{\ast}}(\mathcal{A})\not =0$. Then $\mathcal{A}$ has a computable
copy whose index can be computed effectively from indices for $\mathcal{A}%
_{\omega}$, $R_{\mathcal{A}}$, and~$f_{\mathcal{A}}$.
\end{lem}

The results above, in relativized form, yield the following two theorems,
which we shall use to calculate the complexity of index sets.

\begin{thm}
[Khisamiev]\label{Kh3}

Let $\mathcal{A}$ be a reduced Abelian $p$-group of length $\omega M$, where
$M\in\omega$. Then $\mathcal{A}$ has a computable copy iff for each $k<M:$

\begin{enumerate}
\item the relation $R_{\mathcal{A}}^{k}=\{(r,t):$ $u_{\omega k+r}%
(\mathcal{A})\geq t\}\}$ is $\Sigma_{2k+2}^{0}$, and

\item there is a $\Delta_{2k+1}^{0}$ function $f_{\mathcal{A}}^{k}(r,s)$ such
that for each $n$, the function $f_{\mathcal{A}}^{k}(n,s)$ is nondecreasing
and with limit $r^{\ast}\geq r$ such that $u_{\omega k+r^{\ast}}%
(\mathcal{A})\neq0$.
\end{enumerate}

\noindent Moreover, we can pass effectively from $\Sigma_{2k+2}^{0}$ indices
for the relations $R_{\mathcal{A}}^{k}$, and $\Delta_{2k+1}^{0}$ indices for
appropriate functions $f_{\mathcal{A}}^{k}$ to a computable index for a copy
of $\mathcal{A}$.
\end{thm}

\begin{thm}
[Khisamiev]\label{Kh4}

Let $\mathcal{A}$ be a reduced Abelian $p$-group of length greater than
$\omega M$. Suppose $\mathcal{A}_{\omega M}$ is $\Delta_{2M+1}^{0}$, and for
each $k<M:$

\begin{enumerate}
\item the relation $R_{\mathcal{A}}^{k}=\{(r,t):$ $u_{\omega k+r}%
(\mathcal{A})\geq t\}\}$ is $\Sigma_{2k+2}^{0}$, and

\item there is a $\Delta_{2k+1}^{0}$ function $f_{\mathcal{A}}^{k}(r,s)$ such
that for each $n$, $f_{\mathcal{A}}^{k}(n,s)$ is nondecreasing and with limit
$r^{\ast}\geq r$ such that $u_{\omega k+r^{\ast}}(\mathcal{A})\neq0$.
\end{enumerate}

\noindent Then $\mathcal{A}$ has a computable copy, with index computed
effectively from the $\Delta_{2M+1}^{0}$ index for $\mathcal{A}_{\omega M}$,
$\Sigma_{2k+2}^{0}$ indices for $R_{\mathcal{A}}^{k}$, and $\Delta_{2k+1}^{0}$
indices for appropriate functions $f_{\mathcal{A}}^{k}$.
\end{thm}

\subsection{Index sets of groups of small Ulm length}

In \cite{C2}, it is shown that for the countable reduced Abelian $p$-group of
length $\omega M$ with uniformly infinite Ulm invariants, the index set is
$m$-complete $\Pi_{2M+1}^{0}$. It seemed that other complexities might be
possible for groups of the same length. However, it turned out that the index
set for any group of length $\omega M$ is also $m$-complete $\Pi^{0}_{2m+1}$.
The case $M=1$ of the following result was proved in \cite{CCHM}.

\begin{prop}
\label{apgomega}

Let $K$ be the class of reduced Abelian $p$-groups of length $\omega M$, and
let $\mathcal{A} \in K$. Then $I(\mathcal{A})$ is $m$-complete $\Pi^{0}%
_{2M+1}$ within $K$.
\end{prop}

\proof

Let $\mathcal{A}\in K$. First, we show that $\mathcal{A}$ has a computable
$\Pi_{2M+1}$ Scott sentence. There is a computable $\Pi_{2}$ sentence $\theta$
characterizing the Abelian $p$-groups. Next, there is a computable $\Pi
_{2M+1}$ sentence $\lambda$ characterizing the groups which are reduced and
have length at most $\omega M$. For each $\alpha<\omega M$, we can find a
computable $\Sigma_{2M}$ sentence $\varphi_{\alpha,k}$ saying that $u_{\alpha
}(\mathcal{A})\geq k$. The set of these $\Sigma_{2M}$ sentences true in
$\mathcal{A}$ is $\Sigma_{2M}^{0}$. For each $\varphi_{\alpha,k}$, we can find
a computable $\Pi_{2M}$ sentence equivalent to the negation, and the set of
these sentences true in $\mathcal{A}$ is $\Pi_{2M}^{0}$. We have a computable
$\Pi_{2M+1}$ sentence $\upsilon$ equivalent to the conjunction of the
sentences $\pm\varphi_{\alpha,k}$ true in $\mathcal{A}$. Then we have a
computable $\Pi_{2M+1}$ Scott sentence equivalent to $\theta\ \&\ \lambda
\ \&\ \upsilon$. It follows that $I(\mathcal{A})$ is $\Pi_{2M+1}^{0}$.

For completeness, let $S$ be a $\Pi_{2M+1}^{0}$ set. We will produce a
uniformly computable sequence $(\mathcal{A}_{n})_{n\in\omega}$ of elements of
$K$, such that $n\in S$ if and only if $\mathcal{A}_{n}\cong\mathcal{A}$. We
will specify $\mathcal{A}_{n}$ by giving relations $R_{\mathcal{A}_{n}}^{k}$
and functions $f_{\mathcal{A}_{n}}^{k}$, for $k<M$, as in Theorem \ref{Kh3}.
Since $\mathcal{A}$ is computable, there are relations $R_{\mathcal{A}}^{k}$
and functions $f_{\mathcal{A}}^{k}$, as in the theorem. For $k<M-1$, we let
$R_{\mathcal{A}_{n}}^{k}=R_{\mathcal{A}}^{k}$ and $f_{\mathcal{A}_{n}}%
^{k}=f_{\mathcal{A}}^{k}$.

We will define a single $\Delta_{2M-1}^{0}$ function $g$ to serve as
$f_{\mathcal{A}_{n}}^{M-1}$ for all $n$. We then define a family of relations
$R_{n}$ to serve as $R_{\mathcal{A}_{n}}^{M-1}$, with the feature that if
$(m,k+1)$ is included, so is $(m,k)$. The relations will be uniformly
$\Sigma_{2M}^{0}$. Moreover, if $n\in S$, then we will have $R_{n}%
=R_{\mathcal{A}}^{M-1}$, and if $n\notin S$, then $R_{n}$ will include the
pairs $(m,1)$ for $m=\lim_{s}g(r,s)$, but will omit some of the other pairs
$(m,1)$ in $R_{\mathcal{A}}$. Having determined $g$ and $R_{n}$ for
$n\in\omega$, we will be in a position to apply Theorem \ref{Kh3}.

To get the function $g$, we first define a $\Delta_{2M}^{0}$ sequence
$(k_{n})_{n\in\omega}$, where
\[
k_{0}=\lim\nolimits_{s}f_{\mathcal{A}}^{M-1}(0,s)\text{ and }k_{m+1}%
=\lim\nolimits_{s}f_{\mathcal{A}}^{M-1}(k_{m+1},s)\text{.}%
\]
We will define $g(r,s)$ in such a way that $\lim_{s}g(r,s)=k_{2r+1}$. We have
a $\Delta_{2M-1}^{0}$ approximation to the sequence $(k_{n})_{n\in\omega}$.
Let $k_{m,0}=m$ for all $m$, let $k_{0,s+1}=f_{\mathcal{A}}^{M-1}(0,s+~1)$,
and let $k_{m+1,s+1}$ be the maximum of $f_{\mathcal{A}}^{M-1}(k_{m+1,s+1}%
,s+1)$ and $k_{m+1,s}$. For each $m$, the sequence $k_{m,s}$ is nondecreasing
in $s$ and has limit $k_{m}$. Define $g(r,s)=k_{2r+1,s}$.

The relation $R_{\mathcal{A}}^{M-1}$ is c.e.\ in $\Delta_{2M}^{0}$, and we
have a sequence of finite approximations $R_{\mathcal{A},s}^{M-1}$. The set
$S$ is $\Pi_{1}^{0}$ over $\Delta_{2M}^{0}$. We have a $\Delta_{2M}^{0}$
approximation $(S_{s})_{s\in\omega}$ such that if $n\in S$, then $n\in S_{s}$
for all $s$, and if $n\notin S$, then there exists $s_{0}$ such that for
$s<s_{0}$, $n\in S_{s}$, and for $s\geq s_{0}$, $n\notin S_{s}$. We define
uniformly $\Sigma_{2M}^{0}$ relations $R_{n}$. At stage $0$, we have
$R_{n,0}=\emptyset$. At stage $s+1$, we extend $R_{n,s}$ to $R_{n,s+1}$. We
add the pair $(k_{2s+1},1)$. If $n\notin S_{s}$, this is all, but if $n\in
S_{s}$, we include all pairs $(r,k)$ in $R_{\mathcal{A},s}^{M-1}$. If $n\in
S$, then $R_{n}=R_{\mathcal{A}}^{M-1}$. If $n\notin S$, then $R_{n}$ includes
only finitely many of the pairs $(k_{2s},1)$, while $R_{\mathcal{A}}^{M-1}$
includes all of them.

We apply Theorem \ref{Kh3} as planned to obtain a uniformly computable
sequence of reduced Abelian $p$-groups $(\mathcal{A}_{n})_{n\in\omega}$, all
of length $\omega M$, such that $\mathcal{A}_{n}\cong\mathcal{A}$ iff $n\in S$.

\endproof

So far, we have considered groups of limit length $\omega M$. Now, we consider
groups of successor length. There are several cases.

\begin{prop}
\label{Prop5.7}Let $K$ be the class of reduced Abelian $p$-groups of length
$\omega M+N$, where $N>0$. Suppose $\mathcal{A}$ is a computable member of
$K$, where $\mathcal{A}_{\omega M}$ is finite.

\begin{enumerate}
\item If $\mathcal{A}_{\omega M}$ is minimal for the prescribed length,
$($i.e.,\ it is of type $\mathbb{Z}_{p^{N}})$, then $I(\mathcal{A})$ is
$m$-complete $\Pi_{2M+1}^{0}$ within $K$. $($It is $m$-complete $d$%
-$\Sigma_{2M+1}^{0}$ within the class of groups of length $\leq N$.$)$

\item If $\mathcal{A}_{\omega M}$ is not minimal, then $I(\mathcal{A})$ is
$m$-complete $d$-$\Sigma_{2M+1}^{0}$ within $K$.
\end{enumerate}
\end{prop}

\proof

We use the following lemma.

\begin{lem}
\label{LN}

Let $\mathcal{C}$ be a nontrivial finite Abelian $p$-group of length $N$.

\begin{enumerate}
\item If $\mathcal{C}$ is minimal among groups of length $N$, of type
$\mathbb{Z}_{p^{N}}$, then it has a finitary $\Pi_{1}$ Scott sentence within
groups of length $N$. For any $X$, and any set $S$ that is $\Pi_{1}^{0}(X)$,
there is a uniformly $X$-computable sequence $(\mathcal{C}_{n})_{n\in\omega}$
consisting of groups of length $N$ such that
\[
\mathcal{C}_{n}\cong\mathcal{C}\text{ iff }n\in S\text{.}%
\]

\item If $\mathcal{C}$ is not minimal among groups of length $N$, then it has
a finitary $d$-c.e.\ Scott sentence, and for any $X$, and any set $S$ that is
$d$-c.e.\ relative to $X$, there is a uniformly $X$-computable sequence
$(\mathcal{C}_{n})_{n\in\omega}$ consisting of groups of length $N$ such that
\[
\mathcal{C}_{n}\cong\mathcal{C}\text{ iff }n\in S\text{.}%
\]

\end{enumerate}
\end{lem}

\proof

For 1, we have a finitary $\Pi_{1}$ Scott sentence within $K$ saying that
there are no more than $p^{N}$ elements. The construction will be uniform in
$X$, so, without loss of generality, we assume that $X=\emptyset$. If $S$ is
$\Pi_{1}^{0}$, we have a uniformly computable sequence $(\mathcal{C}%
_{n})_{n\in\omega}$ of Abelian $p$-groups, all of length $N$, such that if
$n\in S$, then $\mathcal{C}_{n}\cong\mathbb{Z}_{p^{N}}$, and if $n\notin S$,
then $\mathcal{C}_{n}\cong\mathbb{Z}_{p^{N}}^{2}$.

\bigskip

For 2, we have a finitary $d$-c.e.\ Scott sentence, as in Section \ref{fssec}.
Suppose $S=S_{1}-S_{2}$, where $S_{1}$ and $S_{2}$ are c.e. We let
$\mathcal{C}^{-}$ be a proper subgroup of $\mathcal{C}$, still of length $N$,
and we let $\mathcal{C}^{+}$ be a proper extension of $\mathcal{C}$, also of
length $N$. We get a uniformly computable sequence $(\mathcal{C}_{n}%
)_{n\in\omega}$ such that $\mathcal{C}_{n}\cong\mathcal{C}^{-}$ if $n\notin
S_{1}$, $\mathcal{C}_{n}\cong\mathcal{C}$ if $n\in S_{1}-S_{2}$, and
$\mathcal{C}_{n}\cong\mathcal{C}^{+}$ if $n\in S_{1}\cap S_{2}$. \ \endproof

Now, we can prove Proposition \ref{Prop5.7}. Let $\mathcal{C}=\mathcal{A}%
_{\omega M}$. For 1, we have a computable $\Pi_{2M+1}$ sentence characterizing
the groups $\mathcal{G}$ such that $\mathcal{G}_{\omega M}\cong\mathcal{C}$
within the class of reduced Abelian $p$-groups of length $\omega M+N$. We have
a computable $\Pi_{2M+1}$ sentence characterizing the Abelian $p$-groups
$\mathcal{G}$ such that for all $\alpha<\omega M$, $u_{\alpha}(\mathcal{G}%
)=u_{\alpha}(\mathcal{A})$. The conjunction, equivalent to a computable
$\Pi_{2M+1}$ sentence, is a Scott sentence for $\mathcal{A}$.

For completeness, let $S$ be $\Pi_{2M+1}^{0}$. Note that $S$ is $\Pi_{1}^{0}$
over $\Delta_{2M+1}^{0}$. By Lemma \ref{LN}, we have a uniformly
$\Delta_{2M+1}^{0}$ sequence $(\mathcal{C}_{n})_{n\in\omega}$ of groups of
length $N$ such that $\mathcal{C}_{n}\cong\mathcal{C}$ iff $n\in S$. Since
$\mathcal{A}$ is computable, we get $\Sigma_{2k+2}^{0}$ relations
$R_{\mathcal{A}}^{k}$, and $\Delta_{2k+1}^{0}$ functions $f^{k}$, as required
in Theorem \ref{Kh4}. We obtain a uniformly computable sequence $(\mathcal{A}%
_{n})_{n\in\omega}$ of groups of length $\omega M+N$ such that $\mathcal{A}%
_{n}\cong\mathcal{A}$ iff $n\in S$.\bigskip

For 2, we have a computable $d$-$\Sigma_{2M+1}$ sentence characterizing the
groups $\mathcal{G}$ such that $\mathcal{G}_{\omega M}\cong\mathcal{C}$. We
have a computable $\Pi_{2M+1}$ sentence characterizing the Abelian $p$-groups
such that for all $\alpha<\omega M$, $u_{\alpha}(\mathcal{G})=u_{\alpha
}(\mathcal{A})$. The conjunction, equivalent to a computable $d$%
-$\Sigma_{2M+1}$ sentence, is a Scott sentence for $\mathcal{A}$.

Toward completeness, let $S$ be a $d$-$\Sigma_{2M+1}^{0}$ set. Then $S$ is
$d$-c.e.\ relative to $\Delta_{2M+1}^{0}$. By Lemma \ref{LN}, there is a
uniformly $\Delta_{2M+1}^{0}$ sequence $(\mathcal{C}_{n})_{n\in\omega}$ of
groups of length $N$ such that $\mathcal{C}_{n}\cong\mathcal{C}$ iff $n\in S$.
As above, since $\mathcal{A}$ is computable, we have relations $R_{\mathcal{A}%
}^{k}$ and functions $f^{k}$ for $k<M$, as required in Theorem~\ref{Kh4}. We
get a uniformly computable sequence $(\mathcal{A}_{n})_{n\in\omega}$ of groups
of length $\omega M+N$ such that $\mathcal{A}_{n}\cong\mathcal{A}$ iff $n\in
S$.

\endproof

We continue with reduced Abelian $p$-groups $\mathcal{A}$ of length $\omega
M+N$, but now we suppose that $\mathcal{A}_{\omega M}$ is infinite. This means
that for some $k<N$, $u_{\omega M+k}(\mathcal{A})=\infty$.

\begin{prop}
\label{apgm2}

Let $K$ be the class of reduced Abelian $p$-groups of length $\omega M+$ $N$.
Let $\mathcal{A}$ be a computable member of $K$. If there is a unique $k<N$
such that $u_{\omega M+k}(\mathcal{A})=\infty$, and for all $m<k$ we have
$u_{\omega M+m}(\mathcal{A})=0$, then $I(\mathcal{A})$ is $m$-complete
$\Pi_{2M+2}^{0}$ within $K$.
\end{prop}

\proof

We use the following lemma.

\begin{lem}
\label{LemmaN}

Suppose $\mathcal{C}$ is a reduced Abelian $p$-group of length $N$, where
there is a unique $k<N$ such that $u_{k}(\mathcal{C})=\infty$, and for all
$m<k$, $u_{m}(\mathcal{C})=0$.

\begin{enumerate}
\item The group $\mathcal{C}$ has a computable $\Pi_{2}$ Scott sentence.

\item For any set $X$, if $S$ is $\Pi_{2}^{0}(X)$, there is a uniformly
$X$-computable sequence $(\mathcal{C}_{n})_{n\in\omega}$ of reduced Abelian
$p$-groups, all of length $N$, such that
\[
\mathcal{C}_{n}\cong\mathcal{C}\text{ iff }n\in S\text{.}%
\]

\end{enumerate}
\end{lem}

\proof

We have $\mathcal{C}\cong H\oplus\mathbb{Z}_{p^{k+1}}^{\infty}$, where $H$ is
a finite direct sum of $\mathbb{Z}_{p^{i+1}}$ for $k<i<N$. For part 1 we will
use a version of Ramsey's Theorem. First, there is a computable $\Pi_{2}$
sentence characterizing the Abelian $p$-groups. There is a computable $\Pi
_{1}$ sentence saying that the length is at most $N$. There is a computable
$\Pi_{2}$ sentence saying that $u_{m}(\mathcal{C})=0$ for all $m<k$. There is
a finitary $d$-$\Sigma_{1}$ sentence characterizing the groups $\mathcal{G}$
of length $N$ such that $\mathcal{G}_{k+1}\cong\mathcal{C}_{k+1}$. Finally,
there is a computable $\Pi_{2}$ sentence saying that for all $r$, there exists
a substructure of type $\mathbb{Z}_{p^{k+1}}^{r}$. The conjunction of these is
equivalent to a computable $\Pi_{2}$ sentence. We show that it is a Scott
sentence for $\mathcal{C}$.

We show that if $\mathcal{G}$ is a model of the proposed Scott sentence, then
$u_{k}(\mathcal{G})=\infty$. To show that $u_{k}(\mathcal{G})\geq m$, consider
the set of statements $z_{1}x_{1}+\ldots+z_{m}x_{m}=h$, where $z_{i}%
\in\mathbb{Z}_{p}$ and $h\in\mathcal{G}_{k+1}$. Say the number of these
statements is $r$. By Ramsey's Theorem (the finite version), there exists $M$
such that
\[
M\rightarrow(2m)_{r}^{k}\text{;}%
\]
that is, for any partition of $k$-sized subsets of a set of size $M$ into $r$
classes, there is a set of size $2m$ that is \textquotedblleft
homogeneous\textquotedblright\ in the sense that all $k$-sized subsets lie in
the same class in the partition (for example, see \cite{hodges}). Take a
substructure of $\mathcal{G}$ of type $\mathbb{Z}_{p^{k+1}}^{M}$, and from
each factor $\mathbb{Z}_{p^{k+1}}$, take an element $b_{i}$ of height $k$ and
order $p$. If there is no $m$-sized subset independent over $\mathcal{G}%
_{k+1}$, then for each $m$-sized subset $\{b_{i_{1}},\ldots,b_{i_{m}}\}$ (with
$i_{1}<\ldots<i_{m}$), one of the $r$ statements above is satisfied. We
partition according to the first such statement. Take a homogeneous set of
size $2m$, with all $m$-tuples satisfying the statement $z_{1}x_{1}%
+\ldots+z_{m}x_{m}=h$. We have disjoint $m$-tuples $c_{1},\ldots,c_{m}$ and
$d_{1},\ldots,d_{m}$ such that
\[
z_{1}c_{1}+\ldots+z_{m}c_{m}=z_{1}d_{1}+\ldots+z_{m}d_{m}\text{.}%
\]
This is impossible. Therefore, we have $u_{k}(\mathcal{G})\geq m$.

\bigskip

For 2, assume, without loss of generality, that $X=\emptyset$. Let $S$ be a
$\Pi_{2}^{0}$ set. We can produce a uniformly computable sequence
$(\mathcal{C}_{n})_{n\in\omega}$ such that if $n\in S$, then $\mathcal{C}%
_{n}\cong H\oplus\mathbb{Z}_{p^{k+1}}^{\infty}$, and otherwise, $\mathcal{C}%
_{n}\cong H\oplus\mathbb{Z}_{p^{k+1}}^{r}$ for some finite $r$. We start with
a copy of $H$. At each stage when we believe that $n\in S$, we add a new
direct summand of type $\mathbb{Z}_{p^{k+1}}$. Otherwise, we add nothing. The
resulting sequence has the desired properties. \ 

\endproof

Using Lemma \ref{LemmaN}, we can prove Proposition \ref{apgm2}. Let
$\mathcal{C}=\mathcal{A}_{\omega M}$. We have a computable $d$-c.e.\ Scott
sentence for $\mathcal{C}$. It follows that there is a computable $d$%
-$\Sigma_{2M+1}$ sentence describing the groups $\mathcal{G}$ with
$\mathcal{G}_{\omega M}\cong\mathcal{C}$. We have a computable $\Pi_{2m+1}$
sentence characterizing the Abelian $p$-groups $\mathcal{G}$ such that for
$\alpha<\omega M$, $u_{\alpha}(\mathcal{G})=u_{\alpha}(\mathcal{A})$. There is
a computable $\Pi_{2M+2}$ sentence equivalent to the conjunction, and this is
a Scott sentence for $\mathcal{A}$.

For completeness, let $\mathcal{C}=\mathcal{A}_{\omega M}$. Let $S$ be
$\Pi_{2M+2}^{0}$. Then $S$ is $\Pi_{2}^{0}$ over $\Delta_{2M+1}^{0}$. By Lemma
\ref{LemmaN}, we get a sequence $(\mathcal{C}_{n})_{n\in\omega}$, uniformly
$\Delta_{2M+1}^{0}$, such that for all $n$, $\mathcal{C}_{n}$ has length $N$,
and $\mathcal{C}_{n}\cong\mathcal{C}$ iff $n\in S$. Now, we apply Theorem
\ref{Kh4}. Since $\mathcal{A}$ is computable, we have $\Sigma_{2k+2}^{0}$
relations $R_{\mathcal{A}}^{k}$, and $\Delta_{2k+1}^{0}$ functions
$f_{\mathcal{A}}^{k}$ for all $k<M$. From these, together with the sequence
$(\mathcal{C}_{n})_{n\in\omega}$, we obtain a uniformly computable sequence
$(\mathcal{A}_{n})_{n\in\omega}$ such that for all $n$, $\mathcal{A}_{n}$ has
length $\omega M+N$, and $\mathcal{A}_{n}\cong\mathcal{A}$ iff $n\in S$.

\endproof

\begin{prop}
\label{apgm3}

Let $K$ be the class of reduced Abelian $p$-groups of length $\omega M + N$
for some $M, N \in\omega$. Let $\mathcal{A} \in K$. If there is a unique $k <
N$ such that $u_{\omega M + k}(\mathcal{A}) = \infty$, and for at least one $m
< k$ we have $0 < u_{\omega M + m}(\mathcal{A}) < \infty$, then $I(\mathcal{A}%
)$ is $m$-complete $d$-$\Sigma^{0}_{2M+2}$ within $K$.
\end{prop}

\proof

We use the following lemma.

\begin{lem}
\label{Lemmak}

Suppose $\mathcal{C}$ has length $N$. Suppose there is a unique $k<N$ such
that $u_{k}(\mathcal{C})=\infty$, and for at least one $m<k$, $0<u_{m}%
(\mathcal{C})<\infty$.

\begin{enumerate}
\item The structure $\mathcal{C}$ has a computable $d$-$\Sigma_{2}$ Scott sentence.

\item For any $X$, if $S$ is $d$-$\Sigma_{2}^{0}(X)$, there is a uniformly
$X$-computable sequence $(\mathcal{C}_{n})_{n\in\omega}$ of reduced Abelian
$p$-groups of length $N$ such that
\[
\mathcal{C}_{n}\cong\mathcal{C}\text{ iff }n\in S\text{.}%
\]

\end{enumerate}
\end{lem}

\proof

For 1, the Scott sentence is the same as in Lemma \ref{LemmaN}, except that we
must specify $u_{m}(\mathcal{C})$, for $m<k$. If $u_{m}(\mathcal{C})\not =0$,
we need a computable $d$-$\Sigma_{2}$ sentence.

\bigskip

For 2, assume, without loss of generality, that $X=\emptyset$. Let $S$ be
$d$-$\Sigma_{2}^{0}$, say $S=S_{1}-S_{2}$, where $S_{i}$ is $\Sigma_{2}^{0}$.
Let $\mathcal{C}^{-}$ have the same Ulm sequence as $\mathcal{C}$ except that
$u_{m}(\mathcal{C}^{-})=0$. Let $\mathcal{C}^{+}$ have the same Ulm sequence
as $\mathcal{C}$ except that $u_{m}(\mathcal{C}^{+})>u_{m}(\mathcal{C})$. We
produce a uniformly computable sequence $(\mathcal{C}_{n})_{n\in\omega}$ of
Abelian $p$-groups of length $N$ such that
\[
\mathcal{C}_{n}\cong\left\{
\begin{array}
[c]{ll}%
\mathcal{C}^{-} & \text{if }n\notin S_{1}\text{,}\\
\mathcal{C} & \text{if }n\in S_{1}-S_{2}\text{,}\\
\mathcal{C}^{+} & \text{if }n\in S_{1}\cap S_{2}\text{.}%
\end{array}
\right.
\]

We start with $H$ and add further direct summands $\mathbb{Z}_{p^{i+1}}$. At
stage $s$, if we believe that $n\notin S_{1}$, then we convert any direct
summands of the form $\mathbb{Z}_{p^{m+1}}$ to the form $\mathbb{Z}_{p^{k+1}}%
$. If we believe $n\in S_{1}-S_{2}$, we make the number of direct summands of
the form $\mathbb{Z}_{p^{m+1}}$ match that in $\mathcal{C}$. If at stage
$s-1$, we had none, then we create new ones. If at stage $s-1$, we had too
many, then we retain those from the greatest stage $t<s$ where we had the
right number (or too few), and convert the extra ones to $\mathbb{Z}_{p^{k+1}%
}$. If we believe that $n\in S_{1}\cap S_{2}$, we make the number of direct
summands of the form $\mathbb{Z}_{p^{m+1}}$ match that in $\mathcal{C}^{+}$.
In any case, we add a new direct summand of the form $\mathbb{Z}_{p^{k+1}}$. \ 

\endproof

We turn to the proof of Proposition \ref{apgm3}. Let $\mathcal{C}%
=\mathcal{A}_{\omega M}$. By Lemma \ref{Lemmak}, $\mathcal{C}$ has a
computable $d$-$\Sigma_{2}$ Scott sentence. From this, we get a computable
$d$-$\Sigma_{2M+2}$ sentence describing the groups $\mathcal{G}$ such that
$\mathcal{G}_{\omega M}\cong\mathcal{C}$. We have a computable $\Pi_{2M+1}$
sentence characterizing Abelian $p$-groups with Ulm invariants matching those
of $\mathcal{A}$ for $\alpha<\omega M$. The conjunction, which is equivalent
to a $d$-$\Sigma_{2M+2}$ sentence, is a Scott sentence for $\mathcal{A}$.

For completeness, note that if $S$ is $d$-$\Sigma_{2M+2}^{0}$, then $S$ is
$d$-$\Sigma_{2}^{0}$ relative to $\Delta_{2M+1}^{0}$. Applying Lemma
\ref{Lemmak}, we get a uniformly $\Delta_{2M+1}^{0}$ sequence $(\mathcal{C}%
_{n})_{n\in\omega}$ of groups of length $N$ such that $\mathcal{C}_{n}%
\cong\mathcal{C}$ iff $n\in S$. Now, we apply Theorem \ref{Kh4}, with
$\mathcal{C}_{n}$, together with the $\Sigma_{2k+2}^{0}$ relations
$R_{\mathcal{A}}^{k}$, and the $\Delta_{2k+1}^{0}$ functions $f_{\mathcal{A}%
}^{k}$, for $k<M$. We get a uniformly computable sequence of groups
$(\mathcal{A}_{n})_{n\in\omega}$, all of length $\omega M+N$, such that
$\mathcal{A}_{n}\cong\mathcal{A}$ iff $n\in S$. \ 

\endproof

\begin{prop}
\label{apg4}

Let $K$ be the class of reduced Abelian $p$-groups of length $\omega M+N$ for
some $M,N\in\omega$. Let $\mathcal{A}$ be a computable member of $K$. If there
exist $m<k<N$ such that
\[
u_{\omega M+m}(\mathcal{A})=u_{\omega M+k}(\mathcal{A})=\infty\text{,}%
\]
then $I(\mathcal{A})$ is $m$-complete $\Pi_{2M+3}^{0}$ within $K$.
\end{prop}

\proof

We use the following lemma.

\begin{lem}
\label{Lemmakm}

Let $\mathcal{C}$ be a reduced Abelian $p$-group of length $N$. Suppose $k$ is
greatest such that $u_{k}(\mathcal{C})=\infty$, and there exists $m<k$ such
that $u_{m}(\mathcal{C})=\infty$.

\begin{enumerate}
\item The structure $\mathcal{C}$ has a computable $\Pi_{3}$ Scott sentence.

\item For any $X$, if $S$ is $\Pi_{3}^{0}(X)$, then there is a uniformly
$X$-computable sequence $(\mathcal{C}_{n})_{n\in\omega}$, consisting of groups
of length $N$, such that
\[
\mathcal{C}_{n}\cong\mathcal{C}\text{\ iff }n\in S\text{.}%
\]

\end{enumerate}
\end{lem}

\proof

For 1, we have a finitary $\Pi_{2}$ sentence describing the reduced Abelian
$p$-groups of length $\leq N$. For each $i<N$, if $u_{i}(\mathcal{C})$ is
finite, we have a finitary $d$-$\Sigma_{2}$ sentence specifying the value. If
$u_{i}(\mathcal{C})=\infty$, we have a computable $\Pi_{3}$ sentence saying
this. The conjunction is equivalent to a computable $\Pi_{3}$ sentence, and it
is a Scott sentence for $\mathcal{C}$.\bigskip

For 2, assume, without loss of generality, that $X=\emptyset$. Let $S$ be
$\omega-Cof$. We produce a uniformly computable sequence $(\mathcal{C}%
_{n})_{n\in\omega}$ of groups of length $N$ such that if $n\in S$, then
$\mathcal{C}_{n}\cong\mathcal{C}$, and if $n\notin S$, say $\omega-W_{n}$ has
cardinality $r$, then $u_{m}(\mathcal{C}_{n})=r$. We have a computable
sequence $(F_{s})_{s\in\omega}$ of finite approximations to $\omega-W_{n}$.
Let $F_{0}=\emptyset$. Given $F_{s}$, if there is some $x\in W_{n,s+1}\cap
F_{s}$, then for the least such $x$, we let $F_{s+1}$ consist of all $y<x$ in
$F_{s}$. If there is no such $x$, then take the least $y\notin W_{n,s+1}$ such
that $y\notin F_{s}$, and let $F_{s+1}$ be the result of adding $y$ to $F_{s}%
$. We have $x\in\omega-W_{n}$ iff for all sufficiently large $s$, $x\in F_{s}%
$. Moreover, if $\omega-W_{n}$ is finite, then for infinitely many $s$,
$F_{s}=\omega-W_{n}$.

We may suppose that
\[
\mathcal{C}=H\oplus\mathbb{Z}_{p^{m+1}}^{\infty}\oplus\mathbb{Z}_{p^{k+1}%
}^{\infty}\text{.}%
\]
We construct $\mathcal{C}_{n}$ as follows. We start with a copy of $H$. At
stage $s$, say $F_{s}$ has cardinality $r$, where at stage $s-1$ the
cardinality was $r^{\prime}$. If $r^{\prime}<r$, we add direct summands of the
form $\mathbb{Z}_{p^{m+1}}$ to bring the number up to $r$. If $r^{\prime}>r$,
we keep the direct summands of the form $\mathbb{Z}_{p^{m+1}}$ that we had at
the greatest stage $t<s$, where the number was at most $r$, and we give the
remaining ones the form $\mathbb{Z}_{p^{k+1}}$. In any case, we add at least
one new direct summand of the form $\mathbb{Z}_{p^{k+1}}$.

\endproof

Now, we turn to the proof of Proposition \ref{apg4}. Let $\mathcal{C}%
=\mathcal{A}_{\omega M}$. By Lemma~\ref{Lemmakm}, $\mathcal{C}$ has a
computable $\Pi_{3}$ Scott sentence. It follows that there is a computable
$\Pi_{2M+3}$ sentence characterizing the groups $\mathcal{G}$ such that
$\mathcal{G}_{\omega M}\cong\mathcal{C}$. We have a computable $\Pi_{2M+1}$
sentence characterizing the Abelian $p$-groups $\mathcal{G}$ such that for
$\alpha<\omega M$, $u_{\alpha}(\mathcal{G})=u_{\alpha}(\mathcal{A})$. There is
a computable $\Pi_{2M+3}$ sentence equivalent to the conjunction, and this is
a Scott sentence for $\mathcal{A}$.

For completeness, let $S$ be $\Pi_{2M+3}^{0}$. Then $S$ is $\Pi_{3}^{0}$ over
$\Delta_{2M+1}^{0}$. By Lemma~\ref{Lemmakm}, we have a uniformly
$\Delta_{2M+1}^{0}$ sequence $(\mathcal{C}_{n})_{n\in\omega}$ of groups of
length $N$ such that $\mathcal{C}_{n}\cong\mathcal{C}$ iff $n\in S$. Since
$\mathcal{A}$ is computable, we have relations $R_{\mathcal{A}}^{k}$ and
functions $f_{\mathcal{A}}^{k}$, for $k<M$, as required in Theorem \ref{Kh4}.
We get a uniformly computable sequence $(\mathcal{A}_{n})_{n\in\omega}$ of
groups of length $\omega M+N$, such that $\mathcal{A}_{n}\cong\mathcal{A} $
iff $n\in S$.

\endproof

We can now summarize the results for groups $\mathcal{A}$ where $\lambda
(\mathcal{A})<\omega^{2}$.

\begin{thm}
\label{apgmain}

Let $K$ be the class of reduced Abelian $p$-groups of length $\omega M+~N$ for
some $M, N \in\omega$. Let $\mathcal{A} \in K$.

\begin{enumerate}
\item If $\mathcal{A}_{\omega M}$ is minimal for the given length $($of the
form $\mathbb{Z}_{p^{N}})$, then $I(\mathcal{A})$ is $m$-complete $\Pi
_{2M+1}^{0}$ within $K$.

\item If $\mathcal{A}_{\omega M}$ is finite but not minimal for the given
length, then $I(\mathcal{A})$ is $m$-complete $d$-$\Sigma_{2M+1}^{0}$ within
$K$.

\item If there is a unique $k < N$ such that $u_{\omega M + k}(\mathcal{A}) =
\infty$, and for all $m < k$, $u_{\omega M + m}(\mathcal{A}) = 0$, then
$I(\mathcal{A})$ is $m$-complete $\Pi^{0}_{2M+2}$ within $K$.

\item If there is a unique $k < N$ such that $u_{\omega M + k}(\mathcal{A}) =
\infty$ and for some $m < k$ we have $0 < u_{\omega M + m}(\mathcal{A}) <
\infty$, then $I(\mathcal{A})$ is $m$-complete $d$-$\Sigma^{0}_{2M+2}$
within~$K$.

\item If there exist $m<k<N$ such that $u_{\omega M+m}(\mathcal{A})=u_{\omega
M+k}(\mathcal{A})=\infty$, then $I(\mathcal{A})$ is $m$-complete $\Pi
_{2M+3}^{0}$ within $K$.
\end{enumerate}
\end{thm}

\subsection{Groups of greater Ulm length}

Theorem \ref{apgmain} leaves open the possibility, counterintuitive though it
may be, that there is an Abelian $p$-group of length at least $\omega^{2} $
with an arithmetical index set. The following result rules out this possibility.

\begin{thm}
Let $\mathcal{A}$ be a computable reduced Abelian $p$-group of length greater
than $\omega M$. Then for any $\Delta_{2M+1}^{0}$ set $S$, there is a
uniformly computable sequence $(\mathcal{A}_{n})_{n\in\omega}$ such that
\[
\mathcal{A}_{n}\cong\mathcal{A}\text{ iff }n\in S\text{.}%
\]
That is, $I(\mathcal{A})$ is $\Delta_{2M+1}^{0}$-\emph{hard}.
\end{thm}

\proof

Let $\mathcal{C}=\mathcal{A}_{\omega M}$, and let $\mathcal{C}^{\prime}$ be a
finite reduced Abelian $p$-group, not isomorphic to $\mathcal{C}$. Let $S$ be
$\Delta_{2M+1}^{0}$. We have a uniformly $\Delta_{2M+1}^{0}$ sequence
$(\mathcal{C}_{n})_{n\in\omega}$ such that $\mathcal{C}_{n}\cong\mathcal{C}$
if $n\in S$, and $\mathcal{C}_{n}\cong\mathcal{C}^{\prime}$ otherwise. Since
$\mathcal{A}$ is computable, we have relations $R_{\mathcal{A}}^{k}$ and
functions $f_{\mathcal{A}}^{k}$ as in Theorem \ref{Kh4}. Then we get a
uniformly computable sequence $(\mathcal{A}_{n})_{n\in\omega}$ such that
$\mathcal{A}_{n}\cong\mathcal{A}$ iff $n\in S.$

\endproof

The following corollary is immediate.

\begin{cor}
Let $\mathcal{A}$ be an Abelian $p$-group of length at least $\omega^{2}$.
Then $I(\mathcal{A})$ is not arithmetical.
\end{cor}

\section{Models of the original Ehrenfeucht theory}

\label{ehrsec}

An \emph{Ehrenfeucht theory} is a complete theory $T$ having exactly $n$
nonisomorphic countable models for some finite $n>1$. A well-known result of
Vaught \cite{vaught} shows that $n$ cannot equal $2$. Ehrenfeucht gave an
example for $n=3$. Ehrenfeucht told Vaught about his example, and it is
described in \cite{vaught}. The language of the theory has a binary relation
symbol $<$ and constants $c_{n}$ for $n\in\omega$. The axioms say that $<$ is
a dense linear ordering without endpoints, and the constants are strictly
increasing. The theory $T$ has the following three nonisomorphic countable
models. There is the \emph{prime} model, in which there is no upper bound for
the constants. There is the \emph{saturated} model, in which the constants
have an upper bound but no least upper bound. There is the \emph{middle}
model, in which there is a least upper bound for the constants.

\begin{prop}
Let $K$ be the class of models of the original Ehrenfeucht theory $T$. Let
$\mathcal{A}^{1}$ be the prime model, let $\mathcal{A}^{2}$ be the middle
model, and let $\mathcal{A}^{3}$ be the saturated model.

\begin{enumerate}
\item $I(\mathcal{A}^{1})$ is $m$-complete $\Pi_{2}^{0}$ within $K$.

\item $I(\mathcal{A}^{2})$ is $m$-complete $\Sigma_{3}^{0}$ within $K$.

\item $I(\mathcal{A}^{3})$ is $m$-complete $\Pi^{0}_{3}$ within $K$.
\end{enumerate}
\end{prop}

\proof

For 1, first note that there is a computable $\Pi_{2}$ sentence characterizing
the models of $T$ such that
\[
(\forall x)\bigvee\limits_{n\in\omega}\hspace{-0.15in}\bigvee x<c_{n}\text{.}%
\]
This is a Scott sentence for $\mathcal{A}^{1}$. Therefore, $I(\mathcal{A}%
^{1})$ is $\Pi_{2}^{0}$.

Toward completeness, let $S$ be a $\Pi_{2}^{0}$ set. We will build a uniformly
computable sequence $(\mathcal{A}_{n})_{n\in\omega}$ such that
\[
\mathcal{A}_{n}\cong\left\{
\begin{array}
[c]{ll}%
\mathcal{A}^{1} & \text{if }n\in S\text{,}\\
\mathcal{A}^{2} & \text{otherwise.}%
\end{array}
\right.
\]
We have a computable approximation $(S_{s})_{s\in\omega}$ for $S$ such that
\[
n\in S\text{ iff }n\in S_{s}\text{ for infinitely many }s\text{.}%
\]
For fixed $n$, when $n\notin S$, we build the middle model by creating a least
upper bound for the constants we have placed so far and preserving it
until/unless our approximation changes. When $n\in S_{s}$, we destroy the
current least upper bound and place the next constant at the end of the
ordering. If $n$ is in $S$, then the sequence of constants is cofinal, and we
get a copy of the prime model. If $n$ is not in $S$, then for some stage
$s_{0}$, for all $s\geq s_{0}$, we have $n\notin S_{s}$, and we will preserve
the least upper bound created at stage $s_{0}$. Thus, we get a copy of the
middle model.\bigskip

We turn to 2. First, note that there is a computable $\Sigma_{3}$ Scott
sentence for $\mathcal{A}^{2}$, describing a model of $T$ such that
\[
(\exists x)\,[\bigwedge\limits_{n\in\omega}\hspace{-0.15in}\bigwedge
x>c_{n}\ \&\ (\forall y)\,[(\bigwedge\limits_{n\in\omega}\hspace
{-0.15in}\bigwedge y>c_{n})\rightarrow y\geq x]]\text{.}%
\]
It follows that $I(\mathcal{A}^{2})$ is $\Sigma_{3}^{0}$.

Toward completeness, let $S$ be a $\Sigma_{3}^{0}$ set. We build a uniformly
computable sequence $(\mathcal{A}_{n})_{n\in\omega}$ such that
\[
\mathcal{A}_{n}\cong\left\{
\begin{array}
[c]{ll}%
\mathcal{A}^{2} & \text{if }n\in S\text{,}\\
\mathcal{A}^{3} & \text{otherwise.}%
\end{array}
\right.
\]
Note that $S$ is $\Sigma_{2}^{0}$ over $\Delta_{2}^{0}$. We have a $\Delta
_{2}^{0}$ approximation $(S_{k})_{k\in\omega}$ such that $n\in S$ iff for all
sufficiently large $k$, $n\in S_{k}$. Fix $n$. Then there is a $\Delta_{2}%
^{0}$ sequence of instructions $(i_{k})_{k\in\omega}$. We start with an upper
bound for the constants. If $n\notin S_{k}$, then $i_{k}$ says to destroy the
current least upper bound for the constants, moving left, closer to the
constants. If $n\in S_{k}$, then $i_{k}$ says to preserve the current least
upper bound for the constants.

Now, we build the computable model $\mathcal{A}_{n}$ based on approximations
of the sequence of instructions. There are mistakes of two kinds. We may
wrongly guess that $i_{k}$ said to preserve the current least upper bound for
the constants. The result is a delay. We may wrongly guess that $i_{k}$ said
to destroy the current least upper bound for the constants. Having introduced
a new upper bound to the left of this one, we correct our mistake by putting
the next constant to the right of any added elements, so as to preserve the
upper bound as in the instruction.

Again, if $n$ is in $S$, $\Delta_{2}^{0}$ will eventually think so, and we
will eventually preserve a particular least upper bound for the constants,
building the middle model. Otherwise, infinitely often we will create a new
upper bound for the constants, moving to the left, closer to the constants.
The result is the saturated model.

\bigskip

Finally, we turn to 3. We have a $\Pi_{3}^{0}$ Scott sentence for
$\mathcal{A}^{3}$, describing a model of $T$ such that
\[
(\exists x)\,[\bigwedge\limits_{n\in\omega}\hspace{-0.15in}\bigwedge
x>c_{n}]\ \&\ (\forall y)\,[\bigwedge\limits_{n\in\omega}\hspace
{-0.15in}\bigwedge y>c_{n}\Longrightarrow\exists z[\bigwedge\limits_{n\in
\omega}\hspace{-0.15in}\bigwedge z>c_{n}\ \&\ z<y]]\text{.}%
\]
It follows that $I(\mathcal{A}^{3})$ is $\Pi_{3}^{0}$. For completeness, we
notice that the sequence constructed for Part 2 already serves the purpose.

\endproof


\begin{thebibliography}{99}                                                                                               %


\bibitem {B}E.\ Barker, \textquotedblleft Back and forth relations for reduced
Abelian $p$-groups,\textquotedblright\ \emph{Annals of Pure and Applied Logic}
75 (1995), pp.\ 223--249.

\bibitem {C1}W.\ Calvert, ``The isomorphism problem for classes of computable
fields,'' \emph{Archive for Mathematical Logic} 75 (2004), pp.\ 327--336.

\bibitem {C2}W.\ Calvert, \textquotedblleft The isomorphism problem for
computable Abelian $p$-groups of bounded length,\textquotedblright%
\ \emph{Journal of Symbolic Logic} 70 (2005), pp.\ 331--345.

\bibitem {CCHM}W.\ Calvert, D.\ Cenzer, V.\ Harizanov, and A.\ Morozov,
\textquotedblleft$\Delta_{2}^{0}$ categoricity of Abelian $p$%
-groups,\textquotedblright\ preprint.

\bibitem {CCKM}W.\ Calvert, D.\ Cummins, J.\ F.\ Knight, and S.\ Miller,
\textquotedblleft Comparing classes of finite structures,\textquotedblright%
\ \emph{Algebra and Logic} 43 (2004), pp.\ 374--392.

\bibitem {CMS}B.\ F.\ Csima, A.\ Montalb\'{a}n, and R.\ A.\ Shore,
\textquotedblleft Boolean algebras, Tarski invariants, and index
sets,\textquotedblright\ to appear in the \emph{Notre Dame Journal of Formal
Logic}.

\bibitem {dobritsa}V.\ P.\ Dobritsa, \textquotedblleft Complexity of the index
set of a constructive model,\textquotedblright\ \emph{Algebra and Logic} 22
(1983), pp.\ 269--276.

\bibitem {G-K}S.\ S.\ Goncharov and J.\ F.\ Knight, ``Computable structure and
non-structure theorems,'' \emph{Algebra and Logic}\ 41 (2002), pp.\ 351--373
(English translation).

\bibitem {hodges}W.\ Hodges, \emph{A Shorter Model Theory}, Cambridge
University Press, 1997.

\bibitem {kaplansky}I.\ Kaplansky, \emph{Infinite Abelian Groups}, University
of Michigan Press, Ann Arbor, 1954.

\bibitem {keisler}H.\ J.\ Keisler, \emph{Model Theory for Infinitary Logic},
North-Holland, Amsterdam, 1971.

\bibitem {kh3}N.\ G.\ Khisamiev, \textquotedblleft Constructive Abelian
groups,\textquotedblright\ \emph{Handbook of Recursive Mathematics}
(Yu.\ L.\ Ershov, S.\ S.\ Goncharov, A.\ Nerode, and J.\ B.\ Remmel, editors),
vol.\ 2, North-Holland, Amsterdam, 1998, pp.\ 1177--1231.

\bibitem {kh4}N.\ G.\ Khisamiev, \textquotedblleft Constructive Abelian
$p$-groups,\textquotedblright\emph{Siberian Advances in Mathematics} 2 (1992), pp.\ 68--113.

\bibitem {LS}S.\ Lempp and T.\ Slaman, \textquotedblleft The complexity of the
index sets of $\aleph_{0}$-categorical theories and of Ehrenfeucht
theories,\textquotedblright\ to appear in the \emph{Advances in Logic}
(Proceedings of the North Texas Logic Conference, October 8--10, 2004),
Contemporary Mathematics, American Mathematical Society.

\bibitem {AM}A.\ W.\ Miller, \textquotedblleft On the Borel classification of
the isomorphism class of a countable model,\textquotedblright\ \emph{Notre
Dame Journal of Formal Logic} 24 (1983), pp.\ 22--34.

\bibitem {DM}D.\ E.\ Miller, \textquotedblleft The invariant $\Pi_{\alpha}%
^{0}$ separation principle,\textquotedblright\ \emph{Transactions of the
American Mathematical Society} 242 (1978), pp.\ 185--204.

\bibitem {Nadel}M.\ Nadel, \textquotedblleft Scott sentences and admissible
sets,\textquotedblright\ \emph{Annals of Mathematical Logic} 7 (1974), pp.\ 267--294.

\bibitem {soare}R.\ I.\ Soare, \emph{Recursively Enumerable Sets and Degrees},
Springer-Verlag, Berlin, 1987.

\bibitem {vaught}R.\ L.\ Vaught, \textquotedblleft Denumerable models of
complete theories,\textquotedblright\ \emph{Infinitistic Methods: Proceedings
of the Symposium on Foundations of Mathematics, Warsaw, 1959}, Pergamon Press,
1961, pp.\ 303--231.

\bibitem {White}W.\ White, \textquotedblleft On the complexity of categoricity
in computable structures,\textquotedblright\ \emph{Mathematical Logic
Quarterly} 49 (2003), pp.\ 603--614.

\bibitem {WhiteD}W.\ White, \emph{Characterizations for Computable
Structures}, PhD dissertation, Cornell University, 2000.
\end{thebibliography}
\end{document}